\documentclass[11pt]{amsart} 
\usepackage{verbatim, latexsym, amssymb, amsmath,color}
\usepackage{epsfig}

\def\R{\mathbb R}
\def\RP{\mathbb {RP}}
\def\N{\mathbb N}
\def\Z{\mathbb Z}

\def\r{\rm{rel}}

\def\vol{\mathrm{vol}}
\def\dmn{\mathrm{dmn}}

\def\dist{\mathrm{dist}}

\theoremstyle{remark}

\theoremstyle{definition}

\begin{document}
\title{Weyl law for the volume spectrum}

\author{Yevgeny Liokumovich}\address{Department of Mathematics, Massachusetts Institute of Technology, Cambridge, MA 02139}
\author{Fernando C. Marques}\address{Princeton University \\ Fine Hall \\ Princeton, NJ 08544}
\author{Andr\'{e} Neves}
\address{University of Chicago \\ Department of Mathematics \\ Chicago IL 60637\\ USA /Imperial College London\\ Huxley Building \\ 180 Queen's Gate \\ London SW7 2RH \\ United Kingdom}

\thanks{The article was partly written during the first author's visit to Max Planck Institute for Mathematics at Bonn; he is grateful 
to the Institute for its hospitality. The second author was partly supported by NSF-DMS-1509027 and NSF DMS-1311795. The third author was partly supported by ERC-2011-StG-278940 and EPSRC Programme Grant EP/K00865X/1.}

\email{ylio@mit.edu}
\email{coda@math.princeton.edu}
\email{aneves@uchicago.edu, a.neves@imperial.ac.uk}
\maketitle

\begin{abstract}
Given  $M$ a Riemannian manifold with (possibly empty)
boundary, we show that its volume spectrum $\{\omega_p(M)\}_{p\in\N}$ satisfies a Weyl law that was conjectured  by Gromov.
\end{abstract}

\section{Introduction}

Let $(M,g)$ be a compact Riemannian manifold of dimension $n+1$. It is well known that the eigenvalues of the Laplacian 
have the following min-max characterization:
$$\lambda_p = \inf_{p-\rm{plane}\, Q \subset\, W^{1,2}(M)}\, \sup_{f\in Q-\{0\}} \frac{\int_M {|\nabla f|^2}  dV}{\int_M  {f^2} dV}, \quad p\in\N.$$

In 1911, Weyl (\cite{weyl}) proved an asymptotic formula for the sequence of eigenvalues $\{ \lambda_p \}_{p\in\N}$ that had  a tremendous impact in Mathematics. The celebrated Weyl law  states that $$\lim_{p \rightarrow \infty} \lambda_p p^{-\frac{2}{n+1}}  = a(n)\vol(M)^{-\frac{2}{n+1}},$$
where $a(n)= 4 \pi^2 \vol(B)^{-\frac{2}{n+1}}$ and $B$ is the unit ball in $\R^{n+1}$.

Gromov (\cite{gromov0}, \cite[Section 8]{gromov}, \cite[Section 5.2]{gromov2},
\cite{gromov3})
proposed  a very general framework to study several non-linear analogs of the spectral
problem on $M$. In the case we are interested, the space $W^{1,2}(M)$ is replaced by the space $\mathcal Z_{n}(M; \mathbb{Z}_2)$ of
mod 2 flat $n$-cycles in $M$ (if $M$ has no boundary) and the energy functional is replaced by the volume functional (see Section
\ref{gmt} for precise definitions). Almgren \cite{almgren}  showed there is a weak homotopy equivalence between 
$\mathcal Z_{n}(M; \mathbb{Z}_2)$ and
$\RP^{\infty}$ and thus its cohomology ring has a generator $\bar \lambda\in H^1(\mathcal Z_{n}(M; \mathbb{Z}_2);\Z_2)$. Instead of considering $p$-planes in $W^{1,2}(M)$ one considers $p$-sweepouts, i.e., subsets of $\mathcal Z_{n}(M; \mathbb{Z}_2)$ where the $p$-th cup power $\bar \lambda^p$ does not vanish (see Section \ref{almgren.section} for precise definitions).
The definition of width is similar to the above min-max
characterization of the eigenvalues. The $p$-width of $M$, denoted
by $\omega_p(M)$, is defined as the infimum over all real numbers 
$w$, such that there exists a $p$-sweepout with every element having volume at most $w$ (see Section \ref{widths.introduction} for precise definitions or  \cite{guth} for some motivation).  In the same way that eigenvalues are realized by the energy of eigenfunctions, Almgren--Pitts Min-max Theory says that the widths are realized by the volume of minimal surfaces (with a possibly small singular set). Similar considerations apply when $M$ has boundary but one has to use the space  $\mathcal Z_{n,\r}(M,\partial  M;\Z_2)$ of relative mod 2
flat cycles. 

An insightful idea  of Gromov was to understand that,  using the cohomology structure of  $\mathcal Z_{n}(M; \mathbb{Z}_2)$, many properties of the energy spectrum $\{\lambda_p\}_{p\in\N}$ can be extended to the volume spectrum $\{\omega_p(M)\}_{p\in\N}$. For instance, Gromov  and later Guth (\cite{gromov} or
Guth \cite{guth}) showed  the existence of a constant $C=C(M,g)$ for which
$$
  C^{-1} \vol(M)^{\frac{n}{n+1}} p^{\frac{1}{n+1}} \leq \omega_p(B) \leq 
 C \vol(M)^{\frac{n}{n+1}} p^{\frac{1}{n+1}}
$$
for all $p\in\N$.

The  asymptotic behaviour of the volume spectrum has  been studied in  the  paper by  Guth (\cite{guth}). It has  
been used by Marques and Neves to prove existence of infinitely many
minimal hypersurfaces in manifolds with positive Ricci 
curvature \cite{marques-neves-infinitely}.

Gromov conjectured (\cite[8.4]{gromov}) that the volume spectrum
$\{ \omega_p(M) \}_{p\in\N}$ satisfies a Weyl's asymptotic law.
In this paper we confirm this and show

\subsection{Weyl Law for the Volume Spectrum} \label{main1}
{\em There exists a constant $a(n)>0$ such that, for every compact Riemannian manifold $(M^{n+1},g)$ with
(possibly empty) boundary, we have
 $$\lim_{p \rightarrow \infty} \omega_p(M) p^{-\frac{1}{n+1}}= a(n) \vol(M)^{\frac{n}{n+1}}.$$}
\medskip

After our paper was completed, the Weyl Law for the Volume Spectrum was used in a fundamental way in \cite{irie-marques-neves} to prove the following theorem:
\medskip

{\bf Theorem} (Irie, Marques, Neves, 2017): \textit{Let $M^{n+1}$ be a closed manifold of dimension $(n+1)$, with $3\leq (n+1)\leq 7$. Then for a $C^\infty$-generic Riemannian metric $g$ on $M$, the union of all closed, smooth, embedded minimal hypersurfaces is dense.}
\medskip

This theorem settles the generic case of Yau's Conjecture (\cite{yau1}) about the existence of infinitely many minimal surfaces by proving  that a much stronger property holds true: there are infinitely many closed embedded minimal hypersurfaces intersecting any given ball in $M$.

\subsection{} Let $\Omega$ be a bounded 
open subset of $\mathbb{R}^{n+1}$ with smooth boundary, or more generally,
a Lipschitz domain in the sense of \cite[Definition 2.5]{hofmann}
(this is a weaker regularity condition
for the boundary of $\Omega$, see Section \ref{gmt}). For such domains we have a more general Weyl law,
which applies to the space of cycles of dimension
 $0 < k \leq n$.
Assume in addition that 
\begin{equation}\label{topo.condition}
H_{i}( \Omega,\partial  \Omega;\Z_2)\mbox{ is }\Z_2\mbox{ if }i=n+1\mbox{ and }0\mbox{ if }k<i < n+1.
\end{equation}
Similarly to the case of codimension $1$ we
 can define the $p$-width of dimension $k$,
$\omega_p^k(\Omega)$, to be the min-max quantity corresponding
to the $p$-th cup power of the generator
$\bar\lambda \in H^{n+1-k}(\mathcal Z_{k,\r}( \Omega,\partial  \Omega;\Z_2);\Z_2)$.

\subsection{Weyl Law for Euclidean domains} \label{main2}
{\em There exists a  constant $a(n,k)$ such that,  for every Lipschitz domain $\Omega$ satisfying \eqref{topo.condition}, we have
$$\lim_{p \rightarrow \infty} \omega_p^k(\Omega) p^{-\frac{n+1-k}{n+1}}= a(n,k) \vol(\Omega)^{\frac{k}{n+1}}.$$}
\medskip

In this setting the inequalities proven by Gromov and Guth for the widths \cite{gromov, guth} become, for all $0<k<n$ and $p\in\N$, 
 \begin{equation} \label{Gromov-Guth}
  C(\Omega)^{-1} \vol(\Omega)^{\frac{k}{n+1}} p^{\frac{n+1-k}{n+1}} \leq \omega^k_p(\Omega) \leq 
 C(\Omega) \vol(\Omega)^{\frac{k}{n+1}} p^{\frac{n+1-k}{n+1}}
 \end{equation}
 for some constant $C(\Omega)>0$.  For general Riemannian metrics on $\Omega$
 and $k=n$
 such inequalities do not hold with a constant that depends only on the dimension.
There are examples of metrics $g_i$  on the unit $3$-ball $B^3$ with $\vol(B^3,g_i) = 1 $
and $\omega_1^2(B^3,g_i) \rightarrow \infty$ \cite{papasoglu}.
Upper bounds for the widths can be obtained if the constant $C$ is allowed to depend on the conformal 
class of the manifold (see \cite{gal} and \cite{liokumovich}), similarly 
to the upper bounds obtained by Korevaar for the 
eigenvalues of the Laplacian \cite{korevaar}.

\subsection{Overview of proof}\label{overview}

We start by describing  Weyl law for Lipschitz domains $\Omega$. One of the main tools in the proofs is the 
Lusternick-Schnirelman inequality, which relates the widths of $\Omega$ to the widths of its subsets
(see \cite[8.3]{gromov}, \cite[Section 3]{guth}).
The idea behind this inequality is the following.
Let $U_1$ and $U_2$ be two disjoint domains in $\Omega$ and
let $V_i \subset \mathcal Z_{k,\r}( \Omega,\partial  \Omega;\Z_2)$
denote the set of all cycles whose restriction to $U_i$
has mass strictly less than $\omega^k_1(U_i)$, $i=1,2$. A sweepout of $\Omega$ restricts to a sweepout of $U_i$ and 
therefore the generator $\bar \lambda$
vanishes on $V_i$, $i=1,2.$ If $V_1$ and $V_2$ are open, it follows that
the cup power $\bar\lambda^2$ vanishes on $V_1 \cup V_2$. Thus every $2$-sweepout must have an element that does not lie in $V_1 \cup V_2$ and this means that
$\omega_2^k(\Omega) \geq \omega_1^k(U_1) + \omega_1^k(U_2)$. Iterating this
argument we obtain that for $N$ disjoint domains $\{U_i\}_{i=1}^N$ in $\Omega$ and
any $p \geq \sum _{i=1} ^N p_i$ we have

\begin{equation} \label{ls}
\omega_p^k(\Omega) \geq \sum_{i=1} ^N \omega_{p_i} ^k (U_i).
\end{equation}

We first sketch how to prove the existence of the limit for the standard Euclidean cube
$C$ of volume $1$.
To do this we exploit the fact that $C$ admits self-similar tilings.

Let $\tilde \omega^k_p (C) = \omega^k_p (C) p^{-\frac{n+1-k}{n+1}}$. By the upper bound in \eqref{Gromov-Guth} the sequence
$\{\tilde \omega^k_p (C) \}_{p\in\N}$ is bounded. 
Choose $p_l$ and $q_j$ so that $\tilde \omega^k_{p_l} (C)$ and $\tilde \omega^k_{q_j} (C)$  are, respectively, very close to  $\limsup_{p\to\infty} \tilde \omega_{p}(C)$
and  $\liminf_{p\to\infty} \tilde \omega_{p}(C)$,
and $q_j$ is much bigger than  $p_l$. Consider a maximal packing of $C$ by squares 
$\{C_i^*\}_{i=1}^N$ of volume $p_l/q_j$.
From \eqref{ls} we have $$ \omega^k_{q_j}(C) \geq \sum_{i=1}^N \omega^k_{p_l}(C_i^*)
=N \left(\frac{p_l}{q_j}\right)^{\frac{k}{n+1}}\omega^k_{p_l}(C)$$ and thus
 $$\tilde \omega^k_{q_j}(C) \geq N \frac{p_l}{q_j} \tilde \omega^k_{p_l}(C).$$
The maximal packing condition implies that, as $q_j \rightarrow \infty$, we have $N \frac{p_l}{q_j} \rightarrow 1$ and 
so $\limsup$ and $\liminf$ must coincide. Thus we can denote this limit by a constant $a(n,k)$.

To prove that sequences $\{\tilde \omega^k_p(\Omega) \}$ converge to $a(n,k)$ 
for unit volume domains $\Omega$ in $\mathbb{R}^{n+1}$ we use similar arguments but applied to efficient packings of scaled copies of $C$
 in $\Omega$ and scaled copies of $\Omega$ in $C$.

For a Riemannian manifolds $(M,g)$ a crucial difficulty arises because while one can find efficient packings of cubes (or balls) in $M$, one cannot find efficient packings of copies of $M$ in a cube. The former implies that one can repeat the same type of arguments and show that, with $a(n)=a(n,n)$, 
 $$\liminf_{p \rightarrow \infty} \omega_p(M) p^{-\frac{1}{n+1}}\geq a(n) \vol(M)^{\frac{n}{n+1}}$$
but the latter implies that one needs another idea to prove the reverse inequality.

We do this by subdividing $M$ into $N$ small regions
$\{Q_i\}_{i=1}^N$, which are $(1 + \varepsilon)$-bilipschitz 
diffeomorphic to domains in $\mathbb{R}^{n+1}$. 
Mapping each $Q_i$ to $\mathbb{R}^{n+1}$ by the corresponding
diffeomorphism and connecting the images of $Q_i$ by tubes
of small total volume we obtain a connected domain $\Omega\subset \R^{n+1}$. A $p$-sweepout of $\Omega$ induces, via restriction, $p$-sweepouts of $Q_i$, $i=1,\ldots,N$. The elements in these  restricted sweepouts have boundary in $\partial Q_i$, $i=1,\ldots,N$ and so we cannot add them to make a sweepout of $M$. However, we show that we can turn each restricted $p$-sweepout into  a continuous family of cycles in $Q_i$ if we add to each element some chain that  is entirely contained in  $\partial Q_i$, $i=1,\ldots,N$. Combining these $N$ families of cycles we show that we obtain a $p$-sweepout of $M$  where the  mass of each element,  when compared with the corresponding  element in $\mathcal Z_{n,\r}(\Omega, \partial \Omega;\Z_2)$, has increased  at most by the  volume of $\cup_{i=1}^N \partial Q_i$. As $p \rightarrow \infty$ the volume of $\cup_{i=1}^N \partial Q_i$
is negligible compared to $p^{1/(n+1)}$ yielding the desired upper bound.
This is the only part in the argument where we restrict to the
codimension $1$ case.

We now mention two technical issues that arise in the proof of \eqref{ls}. The first issue is that
the restriction of a cycle in $\mathcal Z_{k,\r}( \Omega,\partial \Omega;\Z_2)$ 
to $U_i$ may not belong to $\mathcal Z_{k,\r}( U_i,\partial  U_i;\Z_2)$
because its boundary might have infinite  mass 
and so we can not conclude that a sweepout of $\Omega$
restricts to a sweepout of $U_i$.
To resolve this we use a perturbation argument 
and slicing theorem from geometric measure theory (see Lemma \ref{restriction.map}). 
The second issue is that mass is not a continuous function
on the space of flat cycles and therefore the subsets 
$V_1$ and $V_2$ need not be open.
To resolve this we follow the strategy in \cite{marques-neves} and consider the finer topology
of the mass norm on the space of relative cycles. We show that restricting to this finer topology does not increase the value of the width
(see Section \ref{approximation.results}) and thus there is no loss of generality.

\subsection{Some questions}We list some open questions. 

The first question is to compute the constants $a(n,k)$.
This is unknown even in the simplest case $n=k=1$.
Potential candidates for the asymptotically optimal families
of sweepouts include nodal sets of eigenfunctions
on the flat disc or the round sphere, or zero sets of harmonic 
polynomials on the flat disc.

The second question is whether the argument for widths
of Riemannian manifolds
can be extended to higher codimension.
Namely, is it true that  for a compact Riemannian manifold
$$\lim_{p \rightarrow \infty} \omega_p^k(M) p^{-\frac{n+1-k}{n+1}}= a(n,k) \vol(M)^{\frac{k}{n+1}}$$
for $k<n$,
where $a(n,k)$ is the constant for the corresponding limit for Euclidean domains? 
That the liminf of the sequence on the left side is greater or equal than the right side is shown in Theorem \ref{weyl.below.thm}. 

In the case of higher codimension, the cohomology ring of the space of relative cycles is richer (see \cite{guth}) and so another question would be to understand the asymptotic limit for the widths associated with Steenrod powers.

\medskip

The paper is organized as follows. In Section \ref{definitions.section}
we give necessary definitions and prove some technical results
that we need for the proof of Lusternick-Schnirelman inequality.
In Section \ref{section.domains} we prove Lusternick-Schnirelman inequality
and Theorem \ref{main2}. In Section 4 we prove Theorem \ref{main1}.

\section{Definitions and setup}\label{definitions.section}

 \subsection{Geometric Measure Theory}\label{gmt}

  Given $m\in\N$, $I^m$ denotes the $m$-dimensional cube $I^m=[0,1]^m$. For each $j\in \N$, $I(1,j)$ denotes the cube complex on $I^1$  whose $1$-cells and $0$-cells (those are sometimes called vertices) are, respectively,  
$$[0,3^{-j}], [3^{-j},2 \cdot 3^{-j}],\ldots,[1-3^{-j}, 1]\quad\mbox{and}\quad [0], [3^{-j}],\ldots,[1-3^{-j}], [1].$$
We denote by $I(m,j)$  the cell complex on $I^m$: 
$$I(m,j)=I(1,j)\otimes\ldots \otimes I(1,j)\quad (\mbox{$m$ times}).$$
Then $\alpha=\alpha_1 \otimes \cdots\otimes \alpha_m$ is a $q$-cell of $I(m,j)$ if and only if $\alpha_i$ is a cell
of $I(1,j)$ for each $i$, and $\sum_{i=1}^m {\rm dim}(\alpha_i) =q$. We often abuse notation by identifying  a $q$-cell $\alpha$ with its support: $\alpha_1 \times \cdots \times \alpha_m \subset I^m$. 

Given  $X$ a cubical subcomplex of $I^m$, the cube complex $X(j)$ is the union of all cells of $I(m,j)$ whose support is contained in some cell of $X$. We use the notation $X(j)_q$ to denote the set of all $q$-cells in $X(j)$.
\medskip

The ambient spaces considered in this paper are  compact Riemannian $(n+1)$-manifolds $(M^{n+1},g)$  with smooth and possibly empty boundary $\partial M$. 
We can always assume $M$ is isometrically embedded in some Euclidean space $\mathbb{R}^Q$. We denote by $B_r(p)$ the Euclidean open ball of radius $r$ centered at $p\in \mathbb{R}^Q$.
 
 When $M$ is a  region of $\R^{n+1}$ we allow for less regularity  and require $M$ to be a compact region with  finite perimeter  with the extra property that for all $p\in \partial  M$ there is a hyperplane $H\subset \R^{n+1}$ containing $p$, $r>0$, and a Lipschitz function $\phi$ defined on $H\cap B_r(p)$ such that, denoting  by $N$ a  normal vector to $H$,
 $$  M\cap B_r(p)=B_{r}(p)\cap\{x+tN: x\in H, t\geq \phi(x)\}.$$
We call these regions  {\em Lipschitz domains} (see \cite[Definition 2.5]{hofmann}).

\medskip

For reasons to be explained in Section \ref{almgren.section}, with $0 \leq k<n+1$ fixed, we also assume  that
\begin{equation}\label{topo2.condition}
H_{i}(M,\partial  M;\Z_2)\mbox{ is }\Z_2\mbox{ if }i=n+1\mbox{ and }0\mbox{ if }k<i < n+1.
\end{equation}
When $k=n$, this amounts to require that $M$ is connected.

 \medskip

 The following definitions can be found in \cite[Section 4.1]{federer}. For every $0\leq k\leq n+1$, $\mathcal R_k(M;\Z_2)$ (or $\mathcal R_k(\partial M;\Z_2)$) denotes the set of $k$-dimensional rectifiable mod $2$ flat chains in $\R^{Q}$  whose support lies in $ M$ (or $\partial  M$). The Radon measure in $M$ associated with $T\in \mathcal R_k(M;\Z_2)$ (or $\mathcal R_k(\partial M;\Z_2)$) is denoted by $||T||$, and its support is denoted by ${\rm supp}(T)$.

 The mass ${\bf M}$ of $T\in \mathcal R_k(M;\Z_2)$ is defined in \cite[p. 358]{federer}.
 With $S,T\in \mathcal R_k(M;\Z_2)$ the {\em flat metric} is given by
\begin{eqnarray*}
\mathcal F(T,S)=\inf\{{\bf M}(Q)+{\bf M}(R): &&T-S=R+\partial Q, \\
&&\hspace{.5cm} R\in \mathcal R_k(M;\Z_2), Q\in \mathcal R_{k+1}(M;\Z_2)\}.
\end{eqnarray*}
Given a Lipschitz map $F:M\rightarrow M$, the push-forward of $T\in \mathcal R_k(M;\Z_2)$ is denoted by  ${F}_{\#}(T)$.

If $k\geq 1$,  ${\bf I}_{k}(M;\Z_2)$ (or ${\bf I}_{k}(\partial  M;\Z_2)$) denotes  those elements of $\mathcal R_k(M;\Z_2)$ whose boundary lies in $\mathcal R_{k-1}(M;\Z_2)$. Finally, we also consider the spaces
 $$\mathcal Z_{k}(M;\Z_2)=\{T\in {\bf I}_{k}(M;\Z_2): \partial T=0\}$$
and
$$\mathcal Z_{k}(M,\partial  M;\Z_2)=\{T\in {\bf I}_{k}(M;\Z_2): {\rm support}(\partial T)\subset \partial  M\}.$$

 \subsection{Relative cycles}
 
 We now describe the space of relative cycles following \cite[Definition 1.20]{almgren}. We say that $T,S \in \mathcal Z_{k}(M,\partial  M;\Z_2)$ are equivalent if  $T-S \in {\bf I}_{k}(\partial  M;\Z_2)$ and the space of such equivalence classes is denoted by $\mathcal Z_{k,\r}(M,\partial  M;\Z_2)$. There is a natural projection map 
$$P:\mathcal Z_{k}(M,\partial  M;\Z_2)\rightarrow\mathcal Z_{k,\r}(M,\partial  M;\Z_2).$$
If $U\subset \R^Q$ is an open set  and $P(T)=P(S)$  then $T\llcorner (U\setminus \partial M)=S\llcorner (U\setminus \partial M)$. 

When $\partial M=\emptyset$ then $\mathcal Z_{k,\r}(M,\partial M;\Z_2)$ is identical to $\mathcal Z_{k}(M;\Z_2)$.
 
The flat metric and the mass norm in the space of relative cycles  are defined, respectively, as 
$$\mathcal F(P(S),P(T))=\inf\{\mathcal F(S+R,T):\, R\in \mathcal {\bf I}_{k}(\partial  M;\Z_2)\}$$
or
$${\bf M}(P(T))=\inf\{{\bf M}(T+R):R\in \mathcal {\bf I}_{k}(\partial  M;\Z_2)\}.$$
These definitions do not depend on the choice of $S, T \in \mathcal Z_{k}(M,\partial  M;\Z_2)$.
The {\em flat topology} on either $\mathcal Z_{k}(M,\partial  M;\Z_2)$ or $\mathcal Z_{k,\r}(M,\partial  M;\Z_2)$  denotes the topology induced by the flat metric. With
the topology of the mass norm, the spaces will be denoted by $\mathcal Z_{k}(M,\partial  M;{\bf M};\Z_2)$ or $\mathcal Z_{k,\r}(M,\partial  M;{\bf M};\Z_2)$.

 To keep notation simple, we denote $P(T)\in \mathcal Z_{k,\r}(M,\partial  M;\Z_2)$  by $T$. 
 
 \subsection{Federer--Fleming Compactness Theorem}\label{FFcompactness}{\em The set
 $$\{T\in \mathcal Z_{k,\r}(M,\partial  M;\Z_2): {\bf M}(T)\leq L\}$$
is compact in the flat topology for all $L>0$.
 }
 \begin{proof}  
 There is an open neighborhood $U\subset \R^Q$ of $\partial M$ and a Lipschitz map $F:U\longrightarrow M$ so that $F(x)=x$ for all $x\in\partial M$. If $\partial M$ is smooth this is obvious  while if $M$ is a Lipschitz domain this follows from \cite[(4.6.7)]{hofmann}. Set $\Lambda_r=\{x\in \R^Q: \dist(x,\partial M)<r\}$ and  choose $r_0$ so that $\Lambda_{r_0}\subset U$.
 
 It suffices to consider a sequence $\{P(T_i)\}_{i\in\N}\subset \mathcal Z_{k,\r}(M,\partial  M;\Z_2)$   such that the mass of $T_i \in \mathcal Z_{k}(M,\partial  M;\Z_2)$ is bounded by $2L$. Slicing \cite[Section 28]{simon} gives us, for all $i\in\N$,  $0<r_i<r_0$ so that 
 \begin{itemize}
 \item  $C_i=T_i\llcorner \Lambda_{r_i}\in {\bf I}_{k}(U;\Z_2)$;
 \item ${\bf M}(\partial C_i-\partial T_i)\leq 2L/{r_0}.$
\end{itemize}
From  the fact that ${F}_{\#}(\partial T_i)=\partial T_i$ we obtain that ${\bf M}(\partial {F}_{\#}(C_i)-\partial T_i)$ is uniformly bounded independently of $i\in\N$.  Moreover   $F_{\#}(C_i)\in  {\bf I}_k(\partial M)$ for all $i\in\N$ and so $S_i=T_i-{F}_{\#}(C_i)$ is a sequence in $\mathcal Z_{k}(M,\partial  M;\Z_2)$ with bounded mass and bounded boundary mass. The result then follows from the classical Federer-Fleming Compactness Theorem.
\end{proof}


The following proposition will also be needed.

 \subsection{Proposition}\label{semicontinuity}{\em  Given $T\in \mathcal Z_{k,\r}(M,\partial  M;\Z_2)$, ${\bf M}(T)=||T||(M\setminus \partial M)$ and the mass is lower semicontinuous with respect to the flat topology in $\mathcal Z_{k,\r}(M,\partial  M;\Z_2)$ meaning that if $U\subset \R^Q$ is an open set and $\{T_i\}_{i\in\N}$ a sequence converging to $T$ in the flat topology then
 $$||T||(U\setminus \partial M)\leq \liminf_{i\to\infty}||T_i||(U\setminus\partial M).$$
 In particular, ${\bf M}(T)\leq \liminf_{i\to\infty}{\bf M}(T_i).$
 }
 \begin{proof} Consider the projection map $F:U\longrightarrow M$ so that $F(x)=x$ for all $x\in\partial M$ and $U\subset \R^Q$ is an open neighborhood of $\partial M$.

 Let $\Omega_r=\{x\in \R^Q: \dist(x,\partial M)>r\}$. The current $T\llcorner \Omega_0$ is also rectifiable \cite[3.8 (3)]{federer-fleming} and hence $S=T-T\llcorner \Omega_0$ is a rectifiable current with support in $\partial M$. From the definition of rectifiable currents, this means there is a sequence of integral Lipschitz chains $S_i\in {\bf I}_k(\R^Q)$ with ${\bf M}(S-S_i)\to 0$  as $i\to\infty$. The currents $\{S_i\}_{i\in\N}$ can be chosen to have support in $U$ and thus, because ${F}_{\#}(S)=S$, we obtain that ${F}_{\#}(S_i)\in  {\bf I}_k(\partial M)$ is such that ${\bf M}(S-{F}_{\#}(S_i))\to 0$  as $i\to\infty$. Therefore ${\bf M}(T-{F}_{\#}(S_i))$ tends to ${\bf M}(T\llcorner \Omega_0)$ and so
 $${\bf M}(P(T))\leq {\bf M}(T\llcorner \Omega_0)=||T||(M\setminus \partial M).$$
The opposite inequality is simple to check.

Without loss of generality we assume that the mass of $\{T_i\}_{i\in\N}$ is uniformly bounded in ${\bf I}_{k}(M;\Z_2)$. From \cite[Proposition 1.16]{almgren} we can choose $\{r_j\}_{j\in\N}$ converging to zero such that, for all $j\in\N$, $T_i\llcorner \Omega_{r_j}\in {\bf I}_{k}(M;\Z_2)$ tends to $T\llcorner \Omega_{r_j}\in {\bf I}_{k}(M;\Z_2)$ in the flat topology  as $i\to\infty$. Thus lower semicontinuity of the mass  for integral currents  implies that 
$$||T||(U\cap \Omega_{r_j})=||T\llcorner \Omega_{r_j}||(U)\leq \liminf_{i\to\infty}||T_i\llcorner \Omega_{r_j}||(U)\leq   \liminf_{i\to\infty}||T_i||(U\setminus\partial M).$$
Making $j\to\infty$ implies the desired result.
 \end{proof}

 \subsection{Almgren Isomorphism}\label{almgren.section} We will be succinct and just describe the main concepts. The reader can see  \cite{almgren, marques-neves-infinitely, marques-neves-cdm, marques-neves-cycles} for some background and explicit constructions.
 
   In \cite{almgren} Almgren constructed,  for each pair of  integers $0\leq k\leq n+1$ and $l\geq 1$, an isomorphism
 $$\Lambda_{l,\r}:\pi_l(\mathcal Z_{k,\r}(M,\partial  M;\Z_2),\{0\}) \rightarrow H_{k+l}(M,\partial  M;\Z_2).$$
 When $\partial M=\emptyset$,  $H_{k+l}(M,\partial M;\Z_2)$  is identical to $H_{k+l}(M;\Z_2)$.

The homological assumptions on $M$ (see Section \ref{gmt}) and  the Almgren Isomorphism Theorem \cite{almgren} imply that the homotopy groups of  $$(\mathcal Z_{k,\r}(M,\partial  M;\Z_2),\{0\})$$   are all  trivial except for the $(n+1-k)$-th one, i.e., $\mathcal Z_{k,\r}(M,\partial  M;\Z_2)$ is an Eilenberg-MacLane space $K(\Z_2,n+1-k)$. Thus  from the Hurewicz Theorem and  Universal Coefficients Theorem one has
$$H^{n+1-k}(\mathcal Z_{k,\r}(M,\partial  M;\Z_2);\Z_2) =\Z_2=\{0,\bar\lambda_k\},$$
$$H^{n+1-k}(\mathcal Z_{k,\r}(M,\partial  M;\Z_2);\Z_2)={\rm Hom}(H_{n+1-k}(\mathcal Z_{k,\r}(M,\partial  M;\Z_2);\Z_2),\Z_2).$$

We now describe the map $$\bar \lambda_k:H_{n+1-k}(\mathcal Z_{k,\r}(M,\partial  M;\Z_2);\Z_2)\rightarrow\Z_2.$$

An element  $\sigma$ in $H_{n+1-k}(\mathcal Z_{k,\r}(M,\partial  M;\Z_2);\Z_2)$ is represented by a continuous map $\Phi:Y\rightarrow  \mathcal Z_{k,\r}(M,\partial  M;\Z_2)$ where $Y$ is a $(n+1-k)$-cubical decomposition of the sphere with  ${\partial Y}=0$.
By Theorem 2.5 of \cite{almgren}, there exists a constant $\nu_{M, \partial M}>0$ such that for all  $l\in\N$ with
\begin{equation}\label{isomorphism.constant}
\mathcal F(\Phi(x),\Phi(y))\leq \nu_{ M,\partial M}\mbox{ for all  adjacent vertices }x,y\in Y(l)_0,
\end{equation}
there exists a chain map $\phi:Y(l)\rightarrow {\bf I}_{\ast}(M;\Z_2)$ of degree $k$ that extends $\Phi$:
\begin{itemize}
\item[(i)] $\phi$ coincides with $\Phi$ on $Y(l)_0$;
\item[(ii)] $\phi(\alpha)\in {\bf I}_{k+p}(M;\Z_2)$ if $\alpha$ is a $p$-cell in $Y(l)_p$; 
\item[(iii)] $\partial \phi(\alpha)=\phi(\partial \alpha)$ in $\mathcal Z_{k+p,\r}(M,\partial  M;\Z_2)$ if $\alpha$ is a $(p+1)$-cell in $Y(l)_{p+1}.$
\end{itemize}
Thus, for every $\alpha\in Y(l)_{n+1-k}$ we obtain $\phi(\alpha)\in {\bf I}_{n+1}(M;\Z_2)$. From $\partial Y=0$ and (iii) we deduce that
\begin{equation}\label{isomorphism.formula}
[\phi]:=\sum_{\alpha\in Y(l)_{n+1-k}} \phi(\alpha)\in \mathcal Z_{n+1,\r}(M,\partial  M;\Z_2).
\end{equation}
From the Constancy Theorem \cite[Theorem 26.27]{simon}, $[\phi]$ is either $ M$ or $0$. In the first case  $\bar \lambda_k (\sigma)=1$ and in the second case $\bar \lambda_k (\sigma)=0$.

  \subsection{Widths}\label{widths.introduction}
  
  
 Let $X$ denote a cubical subcomplex of the $m$-dimensional cube $I^m$. Given $p\in\N$ and an integer $0\leq k< n+1$, a continuous map in the flat topology
  $$\Phi:X\rightarrow  \mathcal Z_{k,\r}(M,\partial  M;\Z_2)$$ 
is called a {\em p-sweepout} if  the p-th cup power of $\lambda_k=\Phi^*(\bar\lambda_k)$ is nonzero in $H^{p(n+1-k)}(X;\Z_2)$. The set of all $p$-sweepouts is closed under homotopies in the flat topology.

We denote by  $\mathcal P^k_p(M)$ the set of all $p$-sweepouts  that are continuous in the flat topology and have no concentration of mass, meaning that (see \cite[Section 3.7]{marques-neves-infinitely}) setting
$${\bf m}(\Phi,r)=  \sup\{||\Phi(x)||(B_r(p)\setminus \partial M) : x\in \dmn(\Phi), p\in  M\}$$
we have 
$\lim_{r\to 0} {\bf m}(\Phi,r)=0.$
The definition above is independent of the representative chosen for the relative cycle $\Phi(x)$. 

 Note that two maps in  $\mathcal P^k_p(M)$ can have different domains. 

\subsection{Lemma}\label{cont.in.mass.no.concentr.}{\em Continuous maps in the mass topology have no concentration of mass.}
\begin{proof} Consider $\bar \Phi:X\longrightarrow \mathcal Z_{k,\r}(M,\partial  M;\Z_2)$ a continuous map in the mass topology.

Choose $\delta>0$. From Proposition \ref{semicontinuity} we have that for all $x\in X$ there is $\Phi(x)\in \mathcal Z_{k}(M,\partial  M;\Z_2)$ with $P(\Phi(x))=\bar \Phi(x)$ and such that $||\Phi(x)||(\partial M)\leq \delta/4$. Thus, for all $x,y\in X$,
$${\bf M}(\Phi(x)-\Phi(y))\leq  {\bf M}(\bar\Phi(x)-\bar\Phi(y))+\delta/2.$$
Given $p\in M$ and $x\in X$, there is $r=r(p,x)>0$  and $U_x\subset X$ an open neighborhood of $x$ so that 
$$||\Phi(y)||(B_r(p))<\delta\quad\mbox{for all } y\in U_x.$$

By compactness, we can select a finite covering $\{B_{r_k}(p_k)\times U_{x_k}\}_{k=1}^N$ of $M\times X$, where $r_k=r(p_k,x_k)/2$. If $\bar r=\min\{r_k\}_{k=1}^N$, then
$$||\Phi(x)||(B_{\bar r}(p))<\delta\mbox{ for all }(p, x)\in M\times X$$
and the result follows.
\end{proof}

 The {\em $p$-width of $ M$ of dimension $k$}  is 
 \begin{equation}\label{rescaled.width}
 \omega^k_p(M)=\inf_{\Phi \in \mathcal P^k_p(M)}\sup\{{\bf M}(\Phi(x)): x\in {\rm dmn}(\Phi)\},
 \end{equation}
where ${\rm dmn}(\Phi)$ is the domain of $\Phi$.  We also set $$\tilde\omega^k_p(M)=p^{-\frac{n+1-k}{n+1}}\omega^k_p(M).$$
 

Gromov and Guth \cite{gromov0, guth} studied the asymptotic behaviour of  $\omega^k_p(M)$ as $p\rightarrow \infty$ and showed:
\subsection{Theorem}\label{gromov-guth}{\em There is  $C_{k,M}>0$  such that
$$
C_{k,M}^{-1} \leq \tilde\omega^k_p(M) \leq C_{k,M}\quad\mbox{ for all }p\in\N.
$$
}

\subsection{Approximation results} \label{approximation.results}

 Given $X$ a cubical  subcomplex, $l\in\N$, an integer $0\leq k< n+1$, and a map $\phi:X(l)_0\rightarrow  \mathcal{Z}_{k,\r}(M,\partial  M;\Z_2)$, we define the {\em fineness} of $\phi$ to be
$${\bf f}(\phi)=\sup\left\{{\bf M}(\phi(x)-\phi(y)): x,y \mbox{ adjacent vertices in } X(l)_0\right\}.$$

\subsection{Theorem}\label{continuous.discrete} {\em Let $\Phi:X\rightarrow \mathcal{Z}_{k,\r}(M,\partial  M;\Z_2)$ be a continuous map in the flat topology  that has no concentration of  mass. 
There exist a sequence of maps
$$\phi_i:X(l_i)_0 \rightarrow \mathcal{Z}_{k,\r}(M,\partial  M;\Z_2), \quad i\in\N$$
with $l_i<l_{i+1}$, and a
sequence of positive numbers $\{\delta_i\}_{i\in\N}$ converging to zero such that  for all $i\in\N$ ${\bf f}(\phi_i)<\delta_i$,
 $$\sup\{\mathcal F(\phi_i(x)-\Phi(x)): x\in X(l_i)_0\}\leq \delta_i,$$
 and
$$\sup\{{\bf M}(\phi_i(x)): x\in X(l_i)_0\}\leq \sup\{{\bf M}(\Phi(x)): x\in X\}+\delta_i.$$
}
  
When $k=n$ and $\partial  M=\emptyset$, the result was proven in  Theorem 13.1 in \cite{marques-neves}.  We leave the proof of Theorem \ref{continuous.discrete} to  Appendix \ref{proofA}.

The purpose of the next theorem is to construct a continuous map in the mass norm out of a discrete map with small fineness.

\subsection{Theorem}\label{discrete.continuous} {\em There exist positive constants $C_0=C_0(M,m)\geq 1$ and $\delta_0=\delta_0(M)$  so that if $Y$ is a cubical subcomplex of $I(m,l)$ and 
$$\phi:Y_0\rightarrow \mathcal{Z}_{k,\r}(M,\partial  M;\Z_2)$$
has ${\bf f}(\phi)<\delta_0$, then there exists   a map 
$$ \Phi:Y\rightarrow \mathcal{Z}_{k,\r}(M,\partial  M;{\bf M};\Z_2)$$
continuous in the mass norm and satisfying  
\begin{itemize}
\item[(i)]$\Phi(x)=\phi(x)$ for all $x\in Y_0$;
\item[(ii)] if $\alpha$ is some $j$-cell in $Y_j$, then $\Phi$ restricted to $\alpha$ depends only on the values of $\phi$ assumed on the vertices of  $\alpha$;
\item[(iii)]
$$\sup\{{\bf M}(\Phi(x)-\Phi(y)): x,y\mbox{ lie in a common cell of } Y\}\leq C_0{\bf f}(\phi).$$
\end{itemize}
}
The map $\Phi$ is called the {\em Almgren extension} of $\phi$. We postpone its proof to Appendix \ref{proofB} because a similar result was proven in Theorem 14.1 of \cite{marques-neves}.
\vskip 0.01in

The proof of Proposition 3.5 of \cite{marques-neves-infinitely} can be extended in a straightforward way to  show that 
\subsection{Proposition}\label{close.implies.homotopic} \textit{Let $Y$ be a cubical subcomplex of  $I(m,l)$.  There exists $\eta=\eta(M,m)>0$ with the following property:}

\textit{ 
If $\Phi_1,\Phi_2:Y \rightarrow \mathcal{Z}_{k,\r}(M,\partial  M;\Z_2)$ are continuous maps in the flat topology such that 
$$\sup\{\mathcal{F}(\Phi_1(y),\Phi_2(y)):y\in Y\} < \eta,$$ then $\Phi_1$ is homotopic to $\Phi_2$ in the flat topology. }
\medskip 

The previous results have the following corollary.

\subsection{Corollary}\label{width.mass}{\em  Given $\Phi\in \mathcal P_p^k(M)$ there is   a
sequence of positive numbers $\{\delta_i\}_{i\in\N}$ converging to zero and a sequence  $\{\Phi_i\}_{i\in\N}\subset \mathcal P_p^k(M)$ of maps continuous in the mass topology such that, for all $i\in\N$,  $\dmn(\Phi)=\dmn(\Phi_i)$  and
$$\sup_{x\in \dmn(\Phi)}\{{\bf M}(\Phi_i(x))\}\leq \sup_{x\in \dmn(\Phi)}\{{\bf M}(\Phi(x))\}+\delta_i.$$
Thus, when computing $\omega^k_p(M)$, we can consider only the maps in $ \mathcal P_p^k(M)$  that are continuous in the mass topology.}

\begin{proof}
From  Theorem \ref{continuous.discrete} and Theorem \ref{discrete.continuous} we obtain a sequence $\{\Phi_i\}_{i\in\N}$ of maps continuous in the mass topology with $\dmn(\Phi)=\dmn(\Phi_i)$ for all $i\in\N$ and a sequence $\{\delta_i\}_{i\in\N}$ tending to zero such that
$$\sup_{x\in \dmn(\Phi)}\{{\bf M}(\Phi_i(x))\}\leq \sup_{x\in \dmn(\Phi)}\{{\bf M}(\Phi(x))\}+\delta_i$$
and
$$\sup_{x\in \dmn(\Phi)}\{\mathcal F(\Phi_i(x)-\Phi(x))\}\leq \delta_i.$$
The corollary follows from Proposition \ref{close.implies.homotopic}.
\end{proof}

\subsection{Restriction of currents}

Let  $R\subset \Omega$ be two Lipschitz domains. Given $T\in \mathcal Z_{k,\r}(\Omega,\partial \Omega;\Z_2)$ it is not necessarily true that $T\llcorner R$ belongs to $\mathcal Z_{k,\r}(R,\partial R;\Z_2)$ because its boundary might have unbounded mass.  Thus the following lemma is needed. The proof is a bit technical and thus could be skipped in a first reading.

\subsection{Lemma}\label{restriction.map}{\em Let $R\subset {\Omega}$ be two  Lipschitz domains.
Consider $0\leq k <n+1$, $p\in\N$, and $\Phi\in\mathcal P^k_{p}(\Omega)$ continuous in the mass topology with $X:=\dmn(\Phi)$.

\begin{itemize}
\item[(1)]
For any $\varepsilon>0$ there is  
$$\hat \Phi:X\rightarrow \mathcal Z_{k,\r}(R,\partial R;\Z_2)$$ continuous in the mass topology  such that  
 $${\bf M}(\hat\Phi(x))\leq (1+\varepsilon)^k{\bf M}(\Phi(x)\llcorner R)\mbox{ for all } x\in X$$
and
$\hat \Phi^{*}\bar \lambda_k=\Phi^*\bar \lambda_k.$

\item[(2)] For any $\bar\varepsilon>0$ and integer $0\leq q\leq p$, the open set
 $$\{x\in X: {\bf M}(\Phi(x)\llcorner {R})<\omega^k_{q}(R)-\bar\varepsilon\}$$
 is contained in an open set $U$  such that  the pull back of $\lambda_k=\Phi^*\bar \lambda_k$ by the inclusion map $\iota: U\rightarrow X$ has its $q$-th cup power vanishing in $H^{q(n+1-k)}(U;\Z_2)$.
 \end{itemize}
 }
\begin{proof}
We start with some discussion regarding the domain $R$. In Theorem 2.7 and Proposition 2.3 of \cite{hofmann} it  is shown the existence of a smooth vector field $Y$ in $\R^{n+1}$ with $|Y|=1$ on $\partial R$ and $\kappa>0$
such that $ Y.\nu\geq \kappa $ a.e. on $\partial R$, where $\nu$ denotes the measure theoretic  inward unit normal of $\partial R$.

Given $t\in \R$ set $$F_t:R \rightarrow \R^{n+1}\quad F_t(x)= x+tY(x).$$
 From \cite[Proposition 4.15 and Theorem 4.19]{hofmann}  we have the existence of  $t_0$ so that,  for all $0<t\leq t_0$, $F_t$ is a bi-Lipschitz map and $R_t:=F_t(R)$  is a    Lipschitz domain contained in the interior of $R$. From  \cite[(4.67)]{hofmann} it also follows the existence of  a Lipschitz function $u$ defined on $\R^{n+1}$ with Lipschitz constant bounded by $C$ such that $R_t=\{x:u(x)\geq t\}$ for all $0\leq t\leq t_0$. Furthermore, we can also assume that 
 $$ |DF_t^{-1}|(x)\leq (1+\varepsilon)\quad\mbox{for all }x\in\Omega.$$

Let $L=\sup_{x\in X}\{{\bf M}(\Phi(x))\}$. For all $l\in\N$ we can find $0<s(l)<t_0$ so that (see \cite[Section 28]{simon})
\begin{itemize}
\item[(a)] $\Phi(x)\llcorner R_{s(l)}\in \mathcal Z_{k,\r}(R_{s(l)},\partial R_{s(l)};\Z_2)$ for all $x\in X(l)_0$;

\medskip
\item[(b)] ${\bf M}(\partial (\Phi(x)\llcorner R_{s(l)}) )\leq CLt_0^{-1}$ for all $x\in X(l)_0$.
\end{itemize}

Set 
$$\phi_l:X(l)_0\rightarrow \mathcal Z_{k,\r}(R,\partial R;\Z_2),\quad \phi_l(x)= (F^{-1}_{s(l)})_{\#}(\Phi(x)\llcorner R_{s(l)}).$$
\medskip
We have for all $x, y \in X(l)_0$
$${\bf M}(\phi_l(x))\leq (1+\varepsilon)^k{\bf M}(\Phi(x)\llcorner R_{s(l)}) \leq (1+\varepsilon)^k {\bf M}(\Phi(x)\llcorner R)$$
and
\begin{equation}\label{mass.diference}
{\bf M}(\phi_l(x)-\phi_l(y))\leq (1+\varepsilon)^k{\bf M}(\Phi(x)-\Phi(y)).
\end{equation}
Moreover
$$\sup_{x\in X}\{{\bf M}(\phi_l(x))+{\bf M}(\partial (\phi_l(x)))\}\leq  (1+\varepsilon)^k(L+CLt_0^{-1})$$
and so we can apply Federer-Fleming Compactness Theorem and a diagonalization argument to conclude the existence of a subsequence $\{\phi_{l_i}\}_{i\in\N}$ such that $\phi_{l_i}(x)$ converges for all $x\in X(l)_0$, $l\in\N$. Using \eqref{mass.diference} we obtain a continuous function in the mass topology
 $$\hat \Phi:X\rightarrow \mathcal Z_{k,\r}(R,\partial R;\Z_2)$$
 such that 
 $${\bf M}(\hat \Phi(x))\leq (1+\varepsilon)^k{\bf M}(\Phi(x)\llcorner R)\mbox{ for all } x\in X.$$
 
For all $i\in\N$ large enough the Almgren extension $\Phi_i$ of $\phi_{l_i}$ is well defined and homotopic in the flat topology to $\hat \Phi$.  Choose such $i\in\N$. To complete the proof of Lemma \ref{restriction.map} (1) we need to show that   $ \Phi_i^{*}\bar \lambda_k=\Phi^*\bar \lambda_k$. 

  An element in $H_{n+1-k}(X;\Z_2)$ is   represented by a continuous map $\tau:Y\rightarrow X$, where $Y$ is a $(n+1-k)$-cubical subcomplex  of $I(m',j)$ and   $\partial Y=0$. 
 
 Choose $i\in\N$ large enough so that for every $x,y$ belonging to a common cell of $X(l_i)$ we have
 \begin{equation}\label{choice.l_i}
{\bf M}(\Phi(x)-\Phi(y))< \min\left\{\nu_{\Omega,\partial \Omega}, \frac{\delta_0}{(1+\varepsilon)^k}, \frac{\nu_{R,\partial R}}{C_0(1+\varepsilon)^k}, \frac{\eta}{5C_0(1+\varepsilon)^k}\right\}
 \end{equation}
 where $\nu_{\Omega,\partial \Omega}, \nu_{R,\partial R}$ are given in \eqref{isomorphism.constant}, $\delta_0=\delta_0(R), C_0=C_0(R,m')$ are given by Theorem \ref{discrete.continuous}, and $\eta=\eta(R,m')$ is given by Proposition \ref{close.implies.homotopic}.  
 
 Choose $t_1$ with $s(l_i)<t_1<t_0$ and such that for all $x\in X(l_i)_0$
   \begin{equation}\label{choice.t_2.2}
  ||\Phi(x)||(R_{s(l_i)}\setminus R_{t_1}) <\frac{\eta}{(1+\varepsilon)^k}
\end{equation} and 
 \begin{equation}\label{choice.t_2.1}
 \mathcal F\left((F^{-1}_{s(l_i)})_{\#}(\Phi(x)\llcorner R_{s(l_i)}),(F^{-1}_{t})_{\#}(\Phi(x)\llcorner R_{s(l_i)})\right)< \frac{\eta}{5}
 \end{equation}
 for all $s(l_i)\leq t\leq t_1$ . The choice in \eqref{choice.t_2.2} is possible because for all $x\in X(l_i)_0$ we have that $||\Phi(x)||(\partial R_{s(l_i)})=0$, and the choice in  \eqref{choice.t_2.1} is possible due to the homotopy formula  \cite[Section 26.22]{simon} . 
 
 Without loss of generality we can assume that for every cell $\alpha \in Y$, $\tau(\alpha)$ is contained in a cell in $X(l_i)_0$.

Set $\Psi=\Phi\circ \tau$ and $ \Psi_i=\Phi_{i}\circ\tau$. These maps represent elements $\sigma$ and $\hat \sigma$ in $H_{n+1-k}(\mathcal Z_{k,\r}( \Omega,\partial  \Omega;\Z_2);\Z_2)$ and $H_{n+1-k}(\mathcal Z_{k,\r}( R,\partial  R;\Z_2);\Z_2)$ respectively.

From \eqref{choice.l_i} we see that the map $\Psi$ satisfies \eqref{isomorphism.constant} (with $l=0$) and so we obtain a chain map of degree $k$ given by 
$\psi:Y\rightarrow {\bf I}_{\ast}( \Omega;\Z_2).$ Using the slicing theory of \cite[Section 28]{simon} we find $s(l_i)<t<t_1$ such that 
\begin{itemize}
\item $\psi(x)\llcorner R_{t}\in \mathcal Z_{k,\r}(R_t,\partial R_t;\Z_2)$ for all $x\in Y_0$;
\item  for every cell $\alpha\in Y_p$ we have $\psi(\alpha)\llcorner R_t\in {\bf I}_{p+k}(R_t;\Z_2)$. 
\end{itemize}
We consider the chain map of degree $k$
$$\bar \psi :Y \rightarrow \mathcal  {\bf I}_{\ast}( R;\Z_2),\quad \bar\psi(\alpha)= (F^{-1}_{t})_{\#}(\psi(\alpha)\llcorner R_{t}).$$
From \eqref{choice.l_i} we see that $\bar \psi_{|Y_0}$ satisfies the conditions of Theorem \ref{discrete.continuous} and we denote by $\bar \Psi$ its Almgren extension. 

From \eqref{choice.l_i} we have that for all $x\in Y$, there are $x'\in X(l_i)_0$ and $x''\in X$ belonging to  a common cell of $X(l_i)$ and such that
$$\mathcal F(\Psi_i(x), \bar \Psi(x))<  \mathcal F\left((F^{-1}_{s(l_i)})_{\#}(\Phi(x')\llcorner R_{s(l_i)}),(F^{-1}_{t})_{\#}(\Phi(x'')\llcorner R_{t})\right)+\frac{2}{5}\eta.$$
Hence we obtain from \eqref{choice.l_i}, \eqref{choice.t_2.2}, and \eqref{choice.t_2.1} that
$$ \mathcal F(\tilde \Psi(x), \bar \Psi(x))<\eta\quad\mbox{for all }x\in Y$$
and Proposition \ref{close.implies.homotopic} implies that $\Psi_i$ and $\bar \Psi$ are homotopic in the flat topology. Hence  they represent the same element $\hat \sigma$ in homology.  

From \eqref{choice.l_i} we see that the map $\bar \Psi$ satisfies \eqref{isomorphism.constant} (with $l=0$). From \eqref{isomorphism.formula} and the definition of $\bar \psi$  we have that $[\psi]\llcorner R= [\bar \psi]$ and so $\bar \lambda_k(\sigma)=\bar \lambda_k(\hat \sigma)$, which is what we wanted to show.

We now prove Lemma \ref{restriction.map} (2). Choose $\varepsilon$ so that 
$$(1+\varepsilon)^k\omega^k_q(R)< \omega^k_q(R)-\bar \varepsilon/2.$$
Considering $\hat\Phi$ given by Lemma \ref{restriction.map} (1) we have that
$$\{x\in X: {\bf M}(\Phi(x)\llcorner R)<\omega^k_q(R)-\varepsilon\}\subset \{x\in X: {\bf M}(\hat\Phi(x))< \omega^k_q(R)-\bar \varepsilon/2\}.$$
Denote the set on the right by $U$ and we assume without loss of generality that its closure $\bar U$ is a cubical complex. If  Lemma \ref{restriction.map} (2) did not hold then $\hat \Phi_{| \bar U} \in \mathcal P^k_{q}(R)$  and this contradicts the definition of $\omega^k_q(R)$.
\end{proof}

\section{Weyl Law for domains} \label{section.domains}

In what follows $C$ denotes the unit cube in $\R^{n+1}$.   Two regions of $\R^{n+1}$ are said to be similar if they differ by an isometry and scaling.
Given a real number, $[x]$ denotes its integer part. Recall the definition of $\tilde \omega_p^k$ in \eqref{rescaled.width}. All domains considered are assumed to satisfy the topological condition \eqref{topo2.condition}.

\subsection{Lusternick-Schnirelman Inequality}\label{LS-inequality}{\em Fix $0\leq k<n+1$ and consider  Lipschitz domains $\Omega_0$, $\{\Omega_i\}_{i=1}^N$, $\{\Omega^*_i\}_{i=1}^N$ such that
\begin{itemize}
\item $\Omega_i$ have unit volume for all $i=0,\ldots, N$;
\item  $\Omega^*_i$ is similar to $\Omega_i$ for all $i=1,\ldots,N$;
\item   $\Omega^*_i\subset {\Omega_0}$ for all $i=1,\ldots,N$ and the interiors of $\{\Omega^*_i\}_{i=1}^N$ are pairwise disjoint.
\end{itemize}
There is a constant $c=c(\Omega_0,k,n)$ such that, with $V=\min\{\vol(\Omega^*_{i})\}_{i=1}^N$ and $p_i=[p\vol(\Omega^*_{i})]$, $i=1,\ldots, N$, we have for all $p\in\N$
$$\tilde \omega^k_p(\Omega_0)\geq \sum_{i=1}^N\vol(\Omega^*_{i})\tilde \omega^k_{p_i}(\Omega_i) -\frac{c}{pV}.$$
}
\begin{proof}
Set $$\bar p=\sum_{i=1}^Np_i=\sum_{i=1}^N[p\vol(\Omega^*_{i})]\leq p\vol(\Omega_0)=p.$$
 Given $\Phi\in\mathcal P^k_p(\Omega_0)$ continuous in the mass topology (with $X= {\rm dmn}(\Phi)$) and $\varepsilon>0$, for each $i=1,\ldots,N$ 
consider the open set $U_i$ given by Lemma \ref{restriction.map} (2) that contains the open set
$$\{x\in X: {\bf M}(\Phi(x)\llcorner \Omega^*_{i})<\omega^k_{p_i}(\Omega^*_i)-\varepsilon/N\}.$$
Thus, denoting by $\iota_i:U_i\rightarrow X$ the inclusion maps,  we have that  $(\iota^*_i \lambda_k)^{p_i}=0$ in $H^{p_i(n+1-k)}(U_i;\Z_2)$ for all $i=1,\ldots,N$, where  $\lambda_k=\Phi^*(\bar \lambda_k)$.

For all $i=1,\ldots,N$, the exact sequence
$$H^{p_i(n+1-k)}(X,U_i;\Z_2)\mathop{\rightarrow}^{j^*} H^{p_i(n+1-k)}(X;\Z_2) \mathop{\rightarrow}^{\iota_{i}^*} H^{p_i(n+1-k)}(U_i ;\Z_2)$$
implies the existence of $\lambda_i\in H^{p_i(n+1-k)}(X,U_i;\Z_2)$ so that $j^*(\lambda_i)=\lambda_k^{p_i}$. Therefore
$$ j^*(\lambda_1)\smile \ldots\smile j^*(\lambda_N)=\lambda_k^{p_1+\ldots+p_N}\neq 0\mbox{ in } H^{\bar p(n+1-k)}(X;\Z_2),$$
because  $\lambda_k^p\neq 0$ and $\bar p\leq p$.

We now claim that $X\neq \cup_{i=1}^NU_i$.  Indeed, if otherwise then
$$H^{\bar p(n+1-k)}(X,U_1\cup \cdots \cup U_N ;\Z_2)=H^{\bar p(n+1-k)}(X,X;\Z_2)=0$$ 
and
from the natural notion of relative cup product (see \cite{hatcher}, p 209)
\begin{multline*}
H^{p_1(n+1-k)}(X,U_1;\Z_2) \smile \cdots \smile H^{p_N(n+1-k)}(X,U_N;\Z_2) \\
\rightarrow H^{\bar p(n+1-k)}(X,U_1\cup \cdots \cup U_N ;\Z_2)=0
\end{multline*}
we see that $\lambda_1 \smile \ldots \smile \lambda_N=0$ which means that
$$
\lambda^{\bar p}=j^*(\lambda_1)\smile \ldots\smile j^*(\lambda_N)=j ^*(\lambda_1 \smile \ldots \smile \lambda_N) =0.
$$
This proves the claim.
\medskip

Thus there is $x\in X\setminus \cup_{i=1}^NU_i$ and so
${\bf M}(\Phi(x))\geq \sum_{i=1}^N \omega^k_{p_i}(\Omega_i^*)-\varepsilon$. Using Corollary \ref{width.mass} and then making $\varepsilon$ tend to zero we obtain
\begin{equation}\label{LS.inequality.noscale}
\omega^k_p(\Omega_0)\geq \sum_{i=1}^N \omega^k_{p_i}(\Omega_i^*).
\end{equation}
As a result there is a  constant $b=b(n,k)$ such that, using Theorem \ref{gromov-guth}, 
\begin{align*}{\tilde \omega^k_p}(\Omega_0)&=p^{-\frac{n+1-k}{n+1}}{\omega^k_p}(\Omega_0)
 \geq p^{-\frac{n+1-k}{n+1}}\sum_{i=1}^N \omega^k_{p_i}(\Omega^*_i)\\
 & = p^{-\frac{n+1-k}{n+1}}\sum_{i=1}^N \vol(\Omega^*_i)^{\frac{k}{n+1}}\omega^k_{p_i}(\Omega_i)\\
&=\sum_{i=1}^N \vol(\Omega^*_i)\left(\frac{p_i}{ p\vol(\Omega^*_i)}\right)^{\frac{n+1-k}{n+1}}\tilde\omega^k_{p_i}(\Omega_i)\\
&\geq \sum_{i=1}^N \vol(\Omega^*_i)\left(1- \frac{1}{ p\vol(\Omega^*_i)}\right)^{\frac{n+1-k}{n+1}}\tilde\omega^k_{p_i}(\Omega_i)\\
&\geq \sum_{i=1}^N \vol(\Omega^*_i)\tilde\omega^k_{p_i}(\Omega_i)-\frac{b}{pV}\sum_{i=1}^N \vol(\Omega^*_i)\tilde\omega^k_{p_i}(\Omega_i)\\
\geq& \sum_{i=1}^N \vol(\Omega^*_i)\tilde\omega^k_{p_i}(\Omega_i)-\frac{bC_k}{pV}\sum_{i=1}^N \vol(\Omega^*_i)\\
\geq & \sum_{i=1}^N \vol(\Omega^*_i)\tilde\omega^k_{p_i}(\Omega_i)-\frac{bC_k}{pV}.
\end{align*}
\end{proof}

We can now prove the main theorem of this section

\subsection{Weyl Law for domains}\label{main.thm.bounded}{\em For all $0\leq k<n+1$, there exists a  constant $a(n,k)$ such that,  for every Lipschitz domain $\Omega$ satisfying \eqref{topo2.condition}, we have
$$\lim_{p \rightarrow \infty} \omega_p^k(\Omega) p^{-\frac{n+1-k}{n+1}}= a(n,k) \vol(\Omega)^{\frac{k}{n+1}}.$$}
\begin{proof}
Without loss of generality we assume that $\Omega$ has unit volume. We start with the following lemma:

\subsection{Lemma}\label{second.lemma.thm}{\em
$\liminf_{p\to\infty}\tilde\omega^k_p(C)= \limsup_{p\to\infty}\tilde\omega^k_p(C).$
}
\begin{proof}
Choose $\{p_l\}_{l\in\N}$, $\{q_j\}_{j\in\N}$ so that
$$\limsup_{p\to\infty}\tilde\omega^k_p(C)=\lim_{l\to\infty}\tilde\omega^k_{p_l}(C)\quad\mbox{and}\quad \liminf_{p\to\infty}\tilde\omega^k_p(C)=\lim_{j\to\infty}\tilde\omega^k_{q_j}(C).$$
With $l$ fixed and for all $j$ large enough so that $\delta_j:= p_l/q_j<1$, consider $N_j$ to be the maximum number of cubes $\{C^*_i\}_{i=1}^{N_j}$,
$\vol(C^*_i) =  \delta_j$ for all $i$, with pairwise disjoint interiors contained in $C$. 
We must have $\delta_j N_j\to 1$ as $j\to\infty$.

From the Lusternick-Schnirelman Inequality \ref{LS-inequality}  we obtain
$$\tilde \omega^k_{q_j}(C)\geq \sum_{i=1}^{N_j}\vol(C^*_{i})\tilde \omega^k_{p_l}(C) + O(p_l^{-1})=\delta_jN_j\tilde \omega^k_{p_l}(C)+O(p_l^{-1}).$$
Making $j\to\infty$ and then $l\to\infty$ we obtain 

$$ \liminf_{p\to\infty}\tilde\omega^k_p(C)\geq \limsup_{p\to\infty}\tilde\omega^k_p(C).$$
\end{proof}

Set $a(n,k)=\lim_{p\to\infty}\tilde\omega^k_p(C)$.

\subsection{Lemma}\label{first.lemma.thm} {\em 
$\liminf_{p\to\infty}\tilde\omega^k_p(\Omega)\geq a(n,k).$
}
\begin{proof}
Given any $\varepsilon>0$, one can find  a family of  cubes $\{C_i^*\}^N_{i=1}$ with pairwise disjoint  interiors contained in $\Omega$, all with the same volume $\delta_i$, and such that $$\sum_{i=1}^N\vol(C_i^*)\geq 1-\varepsilon.$$ From the Lusternick-Schnirelman Inequality \ref{LS-inequality} we obtain
$$\tilde\omega^k_p(\Omega)\geq \sum_{i=1}^N \vol(C^*_i)\tilde\omega^k_{[p\vol(C^*_i)]}(C)-\frac{c}{p\delta_i}$$
and thus making $p\to\infty$ we have
$$\liminf_{p\to\infty}\tilde\omega^k_p(\Omega)\geq (1-\varepsilon)\liminf_{p\to\infty}\tilde\omega^k_p(C)=(1-\varepsilon)a(n,k).$$
The result follows from the arbitrariness of $\varepsilon$.
\end{proof}

\subsection{Lemma}\label{pack.lemma.thm}{\em There are regions $\{\Omega^*_i\}_{i\in\N}$  contained in $C$,  with pairwise disjoint interior, all similar to $\Omega$, and such that for all $\varepsilon>0$ we can choose $N\in\N$ so that
$\sum_{i=1}^N\vol(\Omega^*_i)\geq 1-\varepsilon$. 
}

\begin{proof}
Choose $\Omega_1$ contained in the interior of ${C}$, similar to $\Omega$, and denote its volume by $v$. Set $R_1$ to be the closure $C\setminus \Omega_1$ and find cubes $\{C_{i,1}\}_{i=1}^{Q_1}$ contained in $R_1$ with pairwise disjoint interiors, and such that $\sum_{i=1}^{Q_1}\vol(C_{i,1})\geq \vol(R_1)/2$.  This is possible because $R_1$ is a Lipschitz domain. For all $i=1,\ldots,Q_1$, let $\Omega_{i,1}$ be a region similar to $\Omega_1$, contained in the interior of ${C}_{i,1}$, and with volume $v\vol(C_{i,1})$. 

Next,  set $\Omega_2=\cup_{i=1}^{Q_1}\Omega_{i,1}$ and consider $R_2$ to be the closure of  $C\setminus (\Omega_1\cup \Omega_2).$  Again, find cubes $\{C_{i,2}\}_{i=1}^{Q_2}$ contained in $R_2$, with pairwise disjoint interiors, and such that $\sum_{i=1}^{Q_2}\vol(C_{i,2})\geq \vol(R_2)/2$.  For all $i=1,\ldots,Q_2$, let $\Omega_{i,2}$  be a region similar to $\Omega_1$, contained in the interior of ${C}_{i,2}$, and with volume $v\vol(C_{i,2})$. Define $\Omega_3=\cup_{i=1}^{Q_2}\Omega_{i,2}$ and proceed inductively.

It suffices to check that $\alpha_N:=\sum_{j=1}^N\vol(\Omega_j)$ tends to $1$ as $N$ tends to infinity. Indeed from the construction we have $$\alpha_{N+1}\geq \alpha_N(1-v/2)+v/2\implies \alpha_{N+1}\geq v/2\sum_{j=0}^{N-1}(1-v/2)^j=1-(1-v/2)^N.$$

\end{proof}
\subsection{Lemma}\label{final.lemma.thm}{\em
$a(n,k) \geq \limsup_{p\to\infty}\tilde\omega^k_p(\Omega).$
}
\begin{proof}
Given $\varepsilon>0$, choose $\{q_l\}_{l\in\N}$  so that $$\beta_k:=\limsup_{p\to\infty}\tilde \omega^k_{p}(\Omega)=\lim_{l\to\infty}\tilde \omega^k_{q_l}(\Omega).$$ 

Consider the collection of regions  $\{\Omega^*_i\}_{i=1}^N$ given by the previous lemma. With $l$ fixed and $p$ large, set $\delta_p=q_l/(p\vol (\Omega^*_1))$ and let $Q_p$  be the maximum number of cubes $\{C^*_j\}_{j=1}^{Q_p}$  with pairwise disjoint  interiors contained in $C$ where all have volume $\delta_p$. We  have $\delta_p Q_p$ approaching $1$ as $p \to\infty$.

For each $j=1,\dots, Q_p$ we have regions $\{\Omega_{i,j}\}_{i=1}^N$ inside  $C_j,$ with pairwise disjoint  interiors, all  similar to $\Omega$, and  such that $$\vol(\Omega_{i,j})=\delta_p \vol(\Omega^*_{i})\quad i=1,\ldots,N$$
and thus, with $v=\min\{\vol(\Omega^*_i)\}_{i=1}^N,$
$$\min\{\vol(\Omega_{i,j}):  i=1,\ldots,N, j=1,\dots, Q_p \}=\delta_pv.$$
Set $p_{i}=[p\delta_p \vol(\Omega^*_{i})], i=1,\ldots,N$. From the Lusternick-Schnirelman Inequality \ref{LS-inequality} we have
\begin{align*}
\tilde \omega^k_p(C) & \geq \sum_{j=1}^{Q_p}\vol(\Omega_{1,j})\tilde \omega^k_{p_{1}}(\Omega)+ \sum_{j\geq 1, i\geq 2}\vol(\Omega_{i,j})\tilde \omega^k_{p_{i}}(\Omega)+ O\left(\frac{1}{p\delta_pv}\right)\\
& =Q_{p}\delta_p\left(\vol(\Omega^*_{1})\tilde \omega^k_{q_l}(\Omega)+ \sum_{i=2}^{N}\vol(\Omega^*_{i})\tilde \omega^k_{p_{i}}(\Omega)\right)+O\left(\frac{1}{p\delta_pv}\right)\\
\end{align*}
Making $p$ tend to infinity we obtain 
$$a(n,k)\geq  vol(\Omega^*_1)\tilde \omega^k_{q_l}(\Omega)+  \sum_{i=2}^{N}\vol(\Omega^*_{i})\tilde \omega^k_{[q_l\frac{\vol(\Omega^*_i)}{\vol(\Omega^*_1)}]}(\Omega)+O\left(\frac{\vol(\Omega^*_1)}{q_lv}\right).$$
Making $l\to\infty$ and using Lemma \ref{first.lemma.thm} we obtain that
$$a(n,k)\geq  \vol(\Omega^*_1)\beta_k+a(n,k)\sum_{i=2}^{N}\vol(\Omega^*_{i}).$$
Lemma \ref{pack.lemma.thm} implies then that 
$$ (\vol(\Omega^*_1)+\varepsilon)a(n,k)\geq \vol(\Omega^*_1)\beta_k=\vol(\Omega^*_1)\limsup_{p\to\infty}\tilde \omega^k_{p}(\Omega).$$
The result follows by making $\varepsilon\to 0$.
\end{proof}
The desired result is a consequence of Lemma \ref{first.lemma.thm}, Lemma \ref{second.lemma.thm}, and Lemma \ref{final.lemma.thm}.
\end{proof}

\section{Weyl Law for compact manifolds}

We consider a compact Riemannian $(n+1)$-manifold $(M^{n+1},g)$ isometrically embedded in $\R^Q$  with smooth boundary $\partial M$ and satisfying \eqref{topo2.condition}.  
Recall the definition of $\tilde \omega_p^k(M)$ in \eqref{rescaled.width}, that $C$ denotes the unit cube in $\R^{n+1}$,  and that, for every integer $0\leq k <n+1$, we set $$a(n,k)=\lim_{p\to\infty}\tilde\omega^k_p(C).$$ 
The geodesic ball in $(M,g)$ of radius $r$ centred at $p\in M$ is denoted by $\mathcal B_r(p)$.

\subsection{Theorem}\label{weyl.below.thm}{\em For every integer $0\leq k <n+1$ we have
$$\liminf_{p\to\infty}\omega_p^k(M)p^{-\frac{n+1-k}{n+1}}\geq a(n,k)\vol(M)^{\frac{k}{n+1}}.$$
} 
\begin{proof}
Without loss of generality we assume that $\vol(M)=1$.

Given $\varepsilon>0$ there is $\bar r>0$ so that for every $\mathcal B_r(p)\subset M\setminus \partial M$ with $r\leq \bar r$ and $p\in M$, the Euclidean metric $g_0$ induced on $\mathcal B_r(p)$ via Riemannian normal coordinates centered at $p$ is such that $(1+\varepsilon)^{-1/2}g\leq g_0\leq (1+\varepsilon)^{1/2}g$. Denoting by $|B_r(0)|$ the volume of the Euclidean ball $B_r(0)$ we have
$$(1+\varepsilon)^{-(n+1)}\vol(\mathcal B_r(p))\leq |B_r(0)|\leq (1+\varepsilon)^{n+1}\vol(\mathcal B_r(p))$$
and $\omega^k_p(\mathcal B_r(p))\geq (1+\varepsilon)^{-k}\omega^k_p(B_r(0))$ for all $p\in\N$.

Choose a collection of pairwise disjoint geodesic balls $\mathcal B_i\subset M\setminus \partial M$, $i=1,\ldots, N$, all with radius smaller than $\bar r$, and such that $ \sum_{i=1}^N \vol(\mathcal B_i)\geq (1+\varepsilon)^{-1}$.

Let $B$ denote the unit volume ball in $\R^{n+1}$ and $B_i$ denote a Euclidean ball with the same radius as $\mathcal B_i$, $i=1,\ldots, N$. Reasoning like in the proof of \eqref{LS.inequality.noscale} in Lusternick-Schnirelman Inequality \ref{LS-inequality} we have that
$$\omega^k_p(M)\geq \sum_{i=1}^N \omega^k_{[p\vol(\mathcal B_i)]}(\mathcal B_i)$$ 
and so, with $p_i=[p\vol(\mathcal B_i)]$, $i=1,\ldots,N$,
\begin{align*}
\tilde\omega^k_p (M)&\geq (1+\varepsilon)^{-k}\sum_{i=1}^N |B_i|\left(\frac{p_i}{ p|B_i|}\right)^{\frac{n+1-k}{n+1}}\tilde\omega^k_{p_i}(B)\\
&\geq (1+\varepsilon)^{-k}\sum_{i=1}^N  |B_i|\left(\frac{\vol(\mathcal B_i)}{|B_i|}- \frac{1}{ p|B_i|}\right)^{\frac{n+1-k}{n+1}}\tilde\omega^k_{p_i}(B)\\
&\geq (1+\varepsilon)^{-(n+k+1)}\sum_{i=1}^N  \vol(\mathcal B_i)\left(\frac{\vol(\mathcal B_i)}{|B_i|}- \frac{1}{ p|B_i|}\right)^{\frac{n+1-k}{n+1}}\tilde\omega^k_{p_i}(B).
\end{align*}
Making $p\to\infty$ and using Theorem \ref{main.thm.bounded} we obtain
\begin{multline*}
\liminf_{p\to\infty} \tilde\omega^k_p(M)\geq (1+\varepsilon)^{-2n-2}a(n,k)\sum_{i=1}^N  \vol(\mathcal B_i)
\geq  (1+\varepsilon)^{-2n-3}a(n,k).
\end{multline*}
 The desired result follows by making $\varepsilon$ tend to zero.

\end{proof}

We focus on the case where $k=n$ and  set $a(n)=a(n,n)$. We drop the subscript or superscript $k$ in the notation, which means that $ \omega^n_p(M)$ becomes $\omega_p(M)$, $\mathcal P^n_p(M)$ becomes $\mathcal P_p(M)$ and so on. Condition \eqref{topo2.condition} means that $M$ is connected.

\subsection{Weyl Law for compact manifolds}{\em For every compact Riemannian manifold $(M^{n+1},g)$ with
(possibly empty) boundary, we have
 $$\lim_{p \rightarrow \infty} \omega_p(M) p^{-\frac{1}{n+1}}= a(n) \vol(M)^{\frac{n}{n+1}}.$$}
\begin{proof}

{Recalling the discussion in Section \ref{overview}, we start by decomposing  $M$ into regions that are almost Euclidean (denoted by $\{\mathcal C_i\}_{i=1}^N$)  and then  use those regions to construct a connected region $\Omega$ in Euclidean space.}   

Given $\varepsilon>0$,  consider $\bar r$ so that for all $p\in M$ 
the ball $\mathcal B_{\bar r}(p)$ is $(1+\varepsilon/2)$-bilipschitz diffeomorphic
to some ball of radius $\bar r$ in the closed upper half-space
$\mathbb{R}^n \times \mathbb{R}_+$ with the Euclidean metric.
Choose a covering $\{\mathcal B_i \}_{i= 1}^{\tilde N}$ of $M$ by balls of radius $\bar r$,
so that balls of half the radius still cover $M$.

We now define a collection $\mathcal C=\{\mathcal C_i\}_{i=1}^N$ of domains  with the following properties for all $i=1,\ldots,N$:

\begin{itemize}
\item Each $\mathcal C_i$ is $(1+\varepsilon/2)$-bilipschitz diffeomorphic
to a Lipschitz domain in $\mathbb{R}^{n+1}$ with Euclidean metric;
\item $\mathcal C$ is a covering of $M$;
\item $\mathcal C_i$'s have mutually disjoint interiors.
\end{itemize}

We first define domains $\tilde{\mathcal C_i}$, $i=1,\ldots,\tilde N$, inductively. We set $\tilde {\mathcal C_1} = \mathcal B_1$. For $i>1$ we set $\tilde{\mathcal B_i}$ to be a concentric ball in $\mathcal B_i$ of radius $r_i \in [\frac{\bar r}{2}, \bar r]$,
so that the boundary of $\mathcal B_i$ intersects the boundaries of $\tilde {\mathcal C_1},...,{\tilde {\mathcal C}}_{i-1}$
transversally. 
We define $\tilde {\mathcal C_i}$ to be the closure of  $\tilde{\mathcal B_i} \cap (M \setminus \cup_{j=1}^{i-1} \tilde{ \mathcal C_j})$.
The transversality condition ensures that  $\tilde {\mathcal C_i}$ is a Lipschitz domain for all $i=1,\ldots,\tilde N$. The collection $\mathcal C$ is formed by considering the connected components $\mathcal C_1, \ldots, \mathcal C_N$  of the domains $\tilde{ \mathcal C_1}, \ldots, {\tilde {\mathcal C}}_{\tilde N}$.

For each $i=1,\ldots, N$, let $C_i\subset \R^{n+1}$ be a region $(1+\varepsilon/2)$-bilipschitz diffeomorphic to $\mathcal C_i$. Consider a region $\Omega\subset \R^{n+1}$ that one obtains by 
connecting  the $N$
disjoint regions $C_i\subset\R^{n+1} , i=1,\ldots,N$ 
consecutively by tubes of very small volume. The region $\Omega$ is connected and  a Lipschitz domain.
 Moreover, making the volumes of the connecting tubes sufficiently small
 we obtain 
 
\begin{equation}\label{volume.thm}
\vol(\Omega) \leq (1+\varepsilon)^{n+1} \vol(M).
\end{equation} 

Consider $\Phi\in\mathcal P_p(\Omega)$ continuous in the mass topology with $X=\dmn(\Phi)$.  From Lemma \ref{restriction.map} (i) we obtain, for all $i=1,\ldots,N$, $\Phi_{i}\in \mathcal P_p(C_i)$ with  domain $X$,
\begin{equation}\label{upper.bound.thm}
{\bf M}(\Phi_{i}(x))\leq (1+\varepsilon)^n{\bf M}(\Phi(x)\llcorner C_i)\mbox{ for all } x\in X,
\end{equation}
and $\Phi_{i}^*\bar\lambda=\lambda$, where $\lambda=\Phi^*\bar\lambda$. 

{Next, we  describe  in general terms how to use the maps $\{\Phi_i\}_{i=1}^N$ to  construct a $p$-sweepout   of $M$. The elements   $\Phi_i(x)$ have boundary in $\partial C_i$ and so one can choose $Z_i(x)\in   {\bf I}_{n+1}(C_i;\Z_2)$ so that the cycle $\partial Z_i(x)$  coincides with $\Phi_i(x)$ on the interior of $C_i$. Because the choice of $Z_i(x)$ is not unique  ($C_i+Z_i(x)$ would  have also been a valid choice) it is not always possible to construct a continuous map $x\mapsto \partial Z_i(x)$. Nonetheless, we argue that  a choice of $Z_1$ for a given $x$ induces  choices of $Z_2,\ldots, Z_N$ so that if $\tilde Z_i$  denotes the image of $Z_i$ in $\mathcal C_i$ under the  respective bilipschitz diffeomorphism, then
 $\partial \tilde Z_1+ \ldots +\partial \tilde Z_N$, as a relative cycle of $M$, does not depend on the choice of  $Z_1$.
  Then  we show that the  map $x\mapsto (\partial \tilde Z_1+ \ldots +\partial \tilde Z_N)(x)$ is a $p$-sweepout of $M$  whose elements  have masses comparable  with those of $\Phi$.}
  
For each $i=1,\ldots, N$ set
$$SX_i=\{(x,Z):  x\in X, \Phi_i(x)-\partial Z \in{\bf I}_{n}(\partial C_i;\Z_2)\}\subset X\times {\bf I}_{n+1}(C_i;\Z_2).$$
It is straightforward to see that $SX_i$ does not depend on the choice of the representative for $\Phi_i(x)$ in $\mathcal Z_{n,\r}(C_i,\partial C_i;\Z_2).$ There is a natural projection $\tau_i:SX_i\rightarrow X$, $i=1,\ldots, N$.

\subsection{Lemma} {\em $\tau_i$ is a $2$-cover of $X$ for all $i=1,\ldots,N$.}
\begin{proof}
Fix $i=1,\ldots,N$. 

For every $x\in X$ we have $\Phi_i(x)$ in the connected  component of zero and so from Proposition 1.23 of \cite{almgren} one can find $Z_x\in  {\bf I}_{n+1}(C_i;\Z_2)$ so that  $\Phi_i(x)-\partial Z_x  \in {\bf I}_{n}(\partial C_i;\Z_2)$.  Note that   $(x,C_i+Z_x)$ belongs to $SX_i$ as well.

Given $(x,Z')\in SX_i$, then $\partial (Z'-Z_{x})\in  {\bf I}_{n}(\partial C_i;\Z_2)$ and so we obtain from the Constancy Theorem \cite[page 141]{simon}  that $Z'=Z_x$ or $Z'=C_i+Z_x$. As a result, $\tau_i^{-1}(x)=\{(x,Z_x), (x,C_i+Z_x)\}$.

The unique lifting property holds for $\pi_i$ because  $\mathcal F(Z_x,C_i+Z_x)=|C_i|$  for all $x\in X$ and so for $y$ near $x$ there is a unique  $Z_y$ that is close to $Z_x$ in the flat topology.
\end{proof}

The isomorphism classes of double covers of $X$ are in a bijective correspondence with ${\rm Hom}(\pi_1(X),\Z_2)$, which is homeomorphic to $H^1(X;\Z_2)$.  We claim that, for all $i=1, \ldots, N$,  the element  $\sigma_i\in H^1(X;\Z_2)$ that classifies $SX_i$ is identical to $\lambda$. Indeed given $\gamma:S^1\rightarrow X$ nontrivial in $\pi_1(X)$,  consider a lift to $SX_i$ given by $\theta\mapsto (\gamma(\exp({i\theta})),Z_{\theta})$, $0\leq \theta \leq 2\pi$.  Then $\sigma_i(\gamma)$ is $1$ if $Z_0=C_i-Z_{2\pi}$ and $0$ if $Z_0=Z_{2\pi}$. Thus $\sigma_i(\gamma)$ is non-zero if and only if $\Phi_i\circ\gamma$ is a sweepout.

As a result we obtain that  $SX_1$ is isomorphic to $SX_i$ for all  $i=1,\ldots,N$ and let $F_i:SX_1\rightarrow SX_i$ be the corresponding isomorphism.


For each $i=1,\ldots, N$, there is a natural projection of $SX_i$ into ${\bf I}_{n+1}(C_i;\Z_2)$ that is continuous in the flat topology. Furthermore $C_i$  is bilipschitz diffeomorphic to $\mathcal C_i$ and so, composing the projection map with that diffeomorphism we obtain $\Xi_i:SX_i\rightarrow {\bf I}_{n+1}(\mathcal C_i;\Z_2)$ continuous in the flat topology.

Set
 $$\hat \Psi:SX_1\rightarrow \mathcal Z_{n,\r}(M,\partial M;\Z_2), \quad \hat\Psi(y)=\sum_{i=1}^N\partial(\Xi_i\circ F_i(y)).$$
 The map is continuous in the flat topology. 
 
If $(x,Z)\in SX_1$, then $\Xi_i\circ F_i(x,C_1+Z) =\mathcal C_i+\Xi_i\circ F_i(x,Z) $ for all $i=1,\ldots,N$, and so
\begin{multline*}
\hat\Psi(x,C_1+Z)=\sum_{i=1}^N\partial(\mathcal C_i+\Xi_i\circ F_i(x,Z))=\sum_{i=1}^N\partial \mathcal C_i+\hat\Psi(x,Z) \\
=\partial M+\hat\Psi(x,Z).
\end{multline*}
 Thus $\hat\Psi(x,C_1+Z)=\hat\Psi(x,Z)$ in $ \mathcal Z_{n,\r}(M,\partial M;\Z_2)$, which means that $\hat \Psi$ descends to a  continuous map  in the flat topology $\Psi:X \rightarrow\mathcal  Z_{n,\r}(M,\partial M;\Z_2)$.
 
\subsection{Lemma}\label{sum.thm}{\em For all  $x\in X$ we have
 $${\bf M}(\Psi(x))\leq (1+\varepsilon)^{2n}{\bf M}(\Phi(x))+(1+\varepsilon)^n\sum_{i=1}^N|\partial C_i|.$$
}
\begin{proof}
Choose $(x,Z)\in SX_1$. Then for all $i=1,\ldots,N$,  we have that $F_i(x,Z)=(x,Z_i)\in SX_i$ for some $Z_i\in {\bf I}_{n+1}(C_i;\Z_2)$ and so we deduce from $\partial Z_i-\Phi_i(x)\in {\bf I}_n(\partial C_i,\Z_2)$ and  \eqref{upper.bound.thm} that
$$
{\bf M}(\partial Z_i)  \leq {\bf M}(\Phi_i(x))+|\partial C_i|\leq  (1+\varepsilon)^n{\bf M}(\Phi(x)\llcorner C_i)+|\partial C_i|.
$$
Therefore
$$\sum_{i=1}^N{\bf M}(\partial Z_i)\leq(1+\varepsilon)^n{\bf M}(\Phi(x))+\sum_{i=1}^N|\partial C_i| $$
and the result follows because
${\bf M}(\Psi(x))\leq (1+\varepsilon)^n\sum_{i=1}^N{\bf M}(\partial Z_i).$

\end{proof}

\noindent{\bf Claim :} {\em $\Psi$ is a $p$-sweepout and $\Psi$ has no concentration of mass.}
\begin{proof}
 
Choose $\gamma:S^1\rightarrow X$ nontrivial in $\pi_1(X)$ and denote by $\gamma_1$ its lift to $SX_1$. Then $\gamma_i=F_i\circ \gamma_1$ gives a lift to $SX_i$  for all $i=1,\ldots,N$ and we consider the continuous map in the flat topology
$$B:[0,2\pi]\rightarrow {\bf I}_{n+1}(M;\Z_2),\quad  B(\theta)=\sum_{i=1}^N \Xi_i\circ\gamma_i(\theta).$$
We have  $(\Psi\circ \gamma )(\theta)=\partial B(\theta)$ for all $0\leq \theta\leq 2\pi$.

Hence $\Psi^*\bar\lambda=\lambda$ because, recalling that $\sigma_i=\lambda$ for all $i=1,\ldots,N$,
\begin{align*}
\lambda(\gamma)=0 & \implies \sigma_i(\gamma)=0 \mbox{ for all } i=1,\dots,N \\
& \implies   
 \Xi_i\circ\gamma_i(2\pi)= \Xi_i\circ\gamma_i(0) \mbox{ for all } i=1,\dots,N\\
 &  \implies  B(2\pi)=B(0)
\end{align*}
and
\begin{align*}
\lambda(\gamma)=1 & \implies \sigma_i(\gamma)=1 \mbox{ for all } i=1,\dots,N \\
& \implies   
 \Xi_i\circ\gamma_i(2\pi)=\mathcal C_i+ \Xi_i\circ\gamma_i(0) \mbox{ for all } i=1,\dots,N\\
 &  \implies  B(2\pi)=M+B(0),
\end{align*}
where in the last line we used the fact that $\{\mathcal C_i\}_{i=1}^N$ are pairwise disjoint and cover $M$.

This implies that $\Psi$ is a $p$-sweepout because $\lambda^p\neq 0$.
We leave to the reader to check that $\Psi$ has no concentration of mass.
\end{proof}

From Corollary \ref{width.mass}, Lemma \ref{sum.thm}, and the previous claim we obtain
$$\omega_p(M)\leq (1+\varepsilon)^{2n}\omega_p(\Omega)+ (1+\varepsilon)^n\sum_{i=1}^N|\partial C_i|.$$
Dividing the inequality above by $p^{1/(n+1)}$, making $p\to\infty$, and using  Theorem \ref{main.thm.bounded} we have
$$\limsup_{p\to\infty}\tilde\omega_p(M)\leq a(n)(1+\varepsilon)^{2n}|\Omega|^{n/(n+1)}.$$
Using \eqref{volume.thm}
and making $\varepsilon$ tend to zero in the two inequalities we obtain $$\limsup_{p\to\infty}\tilde\omega_p(M)\leq  a(n)(\vol M)^{n/(n+1)}.$$ This inequality and Theorem \ref{weyl.below.thm} imply the desired result.
\end{proof}

\appendix

\section{}\label{proofA}

\begin{proof}[Proof of Theorem \ref{continuous.discrete}]
 
Set $a(q)=2^{-4(q+2)^2-2}$ where $q\in\N$ is fixed. We use ${\bf B}_{r}^{\mathcal{F}}(T)$ to denote the ball of radius $r$ in the flat topology centred at $T\in \mathcal Z_{k,\r}(M,\partial M;\Z_2)$. Finally $I_0(m,l)$ denotes the cells of $I(m,l)$ whose support lie in $\partial I^m$.

The key step consists  in proving the following lemma below:

Given $T\in \mathcal Z_{k,\r}(M,\partial M;\Z_2)$ with ${\bf M}(T)\leq L$, $l\in\N$, and  $m\leq q+1$, assume there is a sequence
$$\phi_k:I_0(m,l)_0\rightarrow  {\bf B}_{\varepsilon_k}^{\mathcal{F}}(T)\cap\{S:{\bf M}(S)\leq 2L\}$$
with $\varepsilon_k< 1/k$ and ${\bf m}(\phi_k,r)\leq \delta/4$.

\subsection{Lemma}\label{lemma.pitts.singleT} \textit{There exists $N\in \N$, $N \geq l$,  such that for  a subsequence $\{\phi_j\}$ of $\{\phi_k\}_{k\in \mathbb{N}}$ we can find
$$\psi_j:I(1,N)_0\times I_0(m,l)_0\rightarrow  {\bf B}_{\varepsilon_j}^{\mathcal{F}}(T)$$
satisfying
\begin{itemize}
\item[(i)] ${\bf f}(\psi_j) \leq \delta$ if $m=1$ and ${\bf f}(\psi_j)\leq {\bf f}(\phi_j)+\delta$ if $m \neq1$;
\item[(ii)] $\psi_j([0],x)=\phi_j(x)$ and $\psi_j([1],x)=T$ for all $x\in I_0(m,l)_0$;
\item[(iii)] 
\begin{multline*}
\sup\{{\bf M}(\psi_j(y,x)):(y,x)\in I(1,N)_0\times I_0(m,l)_0\}\\
\leq \sup_{x\in  I_0(m,l)_0}\{{\bf M}(\phi_j(x))\}+\frac{\delta}{n+1};
\end{multline*}
\item[(iv)]${\bf m}(\psi_j,r)\leq  2({\bf m}(\phi_j,r)+a(n)\delta).$
\end{itemize}
}
\medskip
Once this result is proven, Theorem \ref{continuous.discrete} follows exactly in the same way that Theorem 13.1 in \cite{marques-neves} followed from \cite[Lemma 13.4]{marques-neves}.

\begin{proof} From Proposition \ref{semicontinuity} we can assume that $||\phi_k(x)||(\partial M)\leq 1/k$ for all $k\in\N$ and $x\in I_0(m,l)_0$. Since the set of varifolds in $\mathcal{V}_k(M)$ with mass bounded above by $2L$ is compact in the weak topology, we can find  a subsequence $\{\phi_j\}$ of $\{\phi_k\}_{k\in \mathbb{N}}$  and  a map
$$V: I_0(m,l)_0\rightarrow \mathcal{V}_k(M)$$
so that 
$$ \lim_{j\to\infty}|\phi_j(x)|=V(x)\mbox{ as varifolds,}$$
for each $x\in  I_0(m,l)_0.$

Note that $\mathcal F(\phi_j(x),T)$ tends to $0$ as $j\to\infty$ (as relative cycles). Thus from Proposition \ref{semicontinuity}  and  since ${\bf m}(\phi_j,r)\leq \delta/4$, we have
\begin{equation}\label{3.7concentration}
||T||(B_r(p)\setminus \partial M)\leq {\bf m}(\phi_j,r)+a(n){\delta}<\frac{\delta}{3}\end{equation}
and
$$ ||V||(B_r(p))\leq {\bf m}(\phi_j,r)+a(n){\delta}<\frac{\delta}{3}$$
for all $j$ sufficiently large,  $p\in M,$ and $x\in I_0(m,l)_0.$

We can choose points $\{p_i\}_{i=1}^v$, and positive real numbers $\{r_i\}_{i=1}^v$, $r_i< r$, so that
$$B_{r_{i_1}}(p_{i_1})\cap B_{r_{i_2}}(p_{i_2})=\emptyset \quad\mbox{if }i_1\neq i_2,$$
and such that
\begin{equation}\label{3.7area.no.boundary}
||T||(\partial B_{r_i}(p_i))\leq ||V(x)||(\partial B_{r_i}(p_i))=0, 
\end{equation}
\begin{equation}\label{3.7area.no.boundary2}
 ||V(x)||(M\setminus \cup_{i=1}^vB_{r_i}(p_i))<\frac{\delta}{3},
\end{equation} 
and
\begin{equation}\label{3.7area.no.concentration}
||T||(B_{r_i}(p_i)\setminus M)\leq ||V(x)||(B_{r_i}(p_i))=\lim_{j\to\infty}||\phi_j(x)||(B_{r_i}(p_i))<\frac{\delta}{3}, 
\end{equation} 
for all $x\in I_0(m,l)_0$ and $i=1,\ldots,v$.  We can assume $v=3^N-1$ for some $N \in \mathbb{N}$ satisfying $N \geq l$.

From  \cite[Proposition 1.23]{almgren}, we get that there exists $Q_j(x)\in {\bf I}_k(M)$, $R_j(x)\in {\bf I}_k(\partial M)$ for all $j$ sufficiently large and  $x\in I_0(m,l)_0$, such that
$$\partial Q_j(x)=\phi_j(x)-T+R_j(x), \quad {\bf M}(Q_j(x))=\mathcal{F}(\phi_j(x)-T).$$
In particular we have ${\bf M}(Q_j(x)) < \varepsilon_j < 1/j$.

For each $i=1,\ldots,v$, consider the distance function $d_i(x)=d(p_i,x)$. Using \cite[Lemma 28.5]{simon}, we  find a decreasing subsequence $\{r^j_i\}$ converging to $r_i$ with $r_i^j<r$ and such that the slices $\langle Q_j(x),d_i,r_i^j\rangle $ are in ${\bf I}_k(M)$ and satisfy
\begin{equation}\label{3.7slice.formula}
\langle Q_j(x),d_i,r_i^j\rangle=\partial (Q_j(x)\llcorner B_{r_i^j}(p_i))-(\phi_j(x)-T+R_j(x))\llcorner B_{r_i^j}(p_i),
\end{equation}
$$\phi_j(x)\llcorner B_{r_i^j}(p_i),\,\, T\llcorner B_{r_i^j}(p_i)\in {\bf I}_k(M), \quad\mbox{and}\quad R_j(x)\llcorner B_{r_i^j}(p_i)\in {\bf I}_k(\partial M)$$
for every $x\in  I_0(m,l)_0$. 
Note that since $\lim_{j\to\infty}{\bf M}(Q_j(x))=0$, by the coarea formula we can choose $\{r^j_i\}$ such that 
\begin{equation}\label{3.7.area.slice}
\sum_{x\in I_0(m,l)_0}\sum_{i=1}^v {\bf M}(\langle Q_j(x),d_i,r_i^j\rangle)\leq a(n){\delta}<\frac{\delta}{2(n+1)}
\end{equation}
for every sufficiently large  $j$.
Furthermore, using  \eqref{3.7area.no.boundary}, \eqref{3.7area.no.boundary2}, \eqref{3.7area.no.concentration}, and Proposition \ref{semicontinuity}, we get that
\begin{equation}\label{3.7area.inequality}
||\phi_j(x)||(B_{r^j_i}(p_i))<\frac{\delta}{3},  \quad ||T||(B_{r^j_i}(p_i)\setminus M)<\frac{\delta}{3},
\end{equation}
\begin{equation}\label{3.7area.inequality.complement}
||\phi_j(x)||(M\setminus \cup_{i=1}^vB_{r_i}(p_i))<\frac{\delta}{3}, \quad ||T||(M\setminus (\cup_{i=1}^vB_{r_i}(p_i)\cup \partial M))<\frac{\delta}{3},
\end{equation}
and
\begin{equation}\label{3.7.lower.inequality}
(||T||-||\phi_j(x)||)(B_{r_i^j}(p_i)\setminus \partial M)\leq \frac{\delta}{2(n+1)v}
\end{equation}
for every sufficiently large $j$, $i=1,\ldots,v,$  and  $x\in I_0(m,l)_0$.

We consider the map given by
\begin{eqnarray*}
\psi_j\left(\left[\frac{i}{3^N}\right],x\right)&=&\phi_j(x)-\sum_{a=1}^i\partial (Q_j(x)\llcorner B_{r^j_a}(p_a))\quad\mbox{ if }0\leq i\leq 3^N-1,\\
 \psi_{j}([1],x)&=&T,
 \end{eqnarray*}
defined on $I(1,N)_0\times I_0(m,l)_0$. 


Note that as relative cycles
$$\psi_j\left(\left[\frac{i}{3^N}\right],x\right)-T=\partial (Q_j(x)\llcorner (M\setminus \cup_{a=1}^iB_{r^j_a}(p_a)),$$
from which it follows that  $\psi_j\left(\left[\frac{i}{3^N}\right],x\right)\in {\bf B}_{\varepsilon_j}^{\mathcal{F}}(T)$. From \eqref{3.7slice.formula},  we also have that as relative cycles
\begin{multline}\label{3.7.formula2}
\psi_j\left(\left[\frac{i}{3^N}\right],x\right)=\phi_j(x)\llcorner (M\setminus \cup_{a=1}^iB_{r^j_a}(p_a))+\sum_{a=1}^iT\llcorner B_{r_a^j}(p_a)\\
-\sum_{a=1}^i\langle Q_j(x),d_a,r_a^j\rangle\llcorner B_{r_a^j}(p_a).
\end{multline}

In what follows the masses of currents are always computed as relative cycles, i..e, using Proposition \ref{semicontinuity}. From  \eqref{3.7.area.slice}, \eqref{3.7area.inequality}, \eqref{3.7area.inequality.complement}, and \eqref{3.7.formula2} we have  that
\begin{eqnarray*}
&&{\bf M}\left(\psi_j\left(\left[\frac{i}{3^N}\right],x\right)-\psi_j\left(\left[\frac{i-1}{3^N}\right],x\right)\right) \\
&&\hspace{2cm}\leq \frac{\delta}{3}+ ||\phi_j(x)||(B_{r_i^j}(p_i))+ ||T||(B_{r_i^j}(p_i)\setminus M)<\delta
\end{eqnarray*}
for $1\leq i\leq v=3^N-1$, and
\begin{multline*}
{\bf M}\left(\psi_j\left(\left[1-\frac{1}{3^N}\right],x\right)-T\right)
\leq ||\phi_j(x)||(M\setminus \cup_{a=1}^vB_{r^j_a}(p_a))\\+ ||T||(M\setminus (\cup_{a=1}^vB_{r^j_a}(p_a)\cup\partial M))+\frac{\delta}{3}<\delta.
\end{multline*}
If ${\bf d}(x,y)=1$, we also have
\begin{eqnarray*}
&&{\bf M}\left(\psi_j\left(\left[\frac{i}{3^N}\right],x\right)-\psi_j\left(\left[\frac{i}{3^N}\right],y\right)\right)\\
&&\hspace{1cm}\leq ||\phi_j(x)-\phi_j(y)||(M\setminus M)+\frac{\delta}{2}\\
&&\hspace{1cm}\leq {\bf f}(\phi_j)+{\delta}.
\end{eqnarray*}
Hence ${\bf f}(\psi_j)\leq {\bf f}(\phi_j)+\delta$.

To prove Lemma \ref{lemma.pitts.singleT}(iii)  we use \eqref{3.7.area.slice}, \eqref{3.7.lower.inequality}, and \eqref{3.7.formula2}, to conclude
\begin{eqnarray*}
{\bf M}\left(\psi_j\left(\left[\frac{i}{3^N}\right],x\right)\right)&\leq& ||\phi_j(x)||(M\setminus (\cup_{a=1}^iB_{r^j_a}(p_a)\cup\partial M))\\
&&\hspace{.5cm}+\sum_{a=1}^i ||T||(B_{r^j_a}(p_a)\setminus \partial M)+\frac{\delta}{2(n+1)}\\ 
&\leq&  ||\phi_j(x)||(M\setminus M)\\
&&\hspace{.5cm}+\sum_{a=1}^i (||T||-||\phi_j(x)||)(B_{r^j_a}(p_a)\setminus M)+\frac{\delta}{2(n+1)}\\
&\leq& ||\phi_j(x)||(M\setminus M)+\frac{\delta}{n+1}= {\bf M}(\phi_j(x))+\frac{\delta}{n+1}.
\end{eqnarray*}
Finally, Lemma \ref{lemma.pitts.singleT}(iv) follows  from \eqref{3.7concentration}, \eqref{3.7.area.slice}, and \eqref{3.7.formula2}:
\begin{multline*}
\left|\left|\psi_j\left(\left[\frac{i}{3^N}\right],x\right)\right|\right|(B_r(p)\setminus \partial M)\leq \\
||\phi_j(x)||(B_r(p)\setminus \partial M)+||T||(B_r(p)\setminus \partial M)+a(n)\delta\\
\leq   2{\bf m}(\phi_k,r)+2a(n){\delta}.
\end{multline*}
\end{proof}

\end{proof}

\section{}\label{proofB}

\begin{proof}[Proof of Theorem \ref{discrete.continuous}]
The analogous result for continuous functions in the flat topology was proven by Almgren in Theorem 6.6 \cite{almgren}. For continuous functions in the mass topology this result was proven in Theorem 14.1 in \cite{marques-neves} when $k=n$ and $\partial  M=\emptyset$ by adapting the proof of \cite[Theorem 6.6]{almgren}.  We now explain which further adaptations need to be made in order to prove Theorem \ref{discrete.continuous}.

The constant $\delta_0$ is chosen so that Theorem 2.5 in \cite{almgren}  can be applied and thus we obtain a chain map (defined in \cite[Definition 2.3]{almgren}) $\phi_{ M}:Y\rightarrow {\bf I}_{\ast}(M;\Z_2)$ of degree $k$ so that $\phi_{ M}=\phi$ on $Y_0$ and 
$${\bf M}(\phi_{ M}(\alpha))\leq 2{\bf f}(\phi)\quad\mbox{for all }\alpha\in Y_p,\, p\geq 1.$$ 
Consider a differentiable triangulation  of $ M$, the deformation map $\mathcal D$ given in \cite[Theorem 4.5]{pitts} that is continuous in the mass topology, and the cutting  function given by \cite[Theorem 5.8]{almgren}. For every $\alpha\in Y_p$, $p\geq 1$, one has now all the necessary ingredients to  consider the function
$$h_{\alpha}:\alpha\rightarrow \mathcal Z_{k}(M,\partial  M;\Z_2)$$
given by \cite[Interpolation Formula 6.3]{almgren} (with $A= M$, $B=\partial  M$, $\phi_A=\phi_{ M}$, and $\phi_B=0$). The projection
$$\tilde h_{\alpha}:\alpha\rightarrow \mathcal Z_{k,\r}(M,\partial  M;{\bf M};\Z_2)$$
is continuous in the mass topology (see \cite[page 297]{almgren}).

Using  the maps $\{\tilde h_{\alpha}\}_{\alpha\in Y}$ and the construction  described in \cite[Section 6.5]{almgren} one obtains the map $\Phi$ satisfying (i) and (ii). Property (iii) also follows because   Theorem 6.6 2 (b) of \cite{almgren} (see also \cite[Lemma 14.4]{marques-neves}) translates into the fact that if $x, y$ lie in a common cell of $Y$ then for some $C=C(M,m)$ we have
$${\bf M}(\Phi(x)-\Phi(y))\leq C\sup\{{\bf M}(\phi_{ M}(\alpha)): \alpha\in Y_p, p\geq 1\}\leq 2C{\bf f}(\phi).$$
\end{proof}

\bibliographystyle{amsbook}

\end{document}